\documentclass[reqno,11pt]{amsart} %maths

\usepackage[top=1.5in,right=1.125in,left=1.125in,bottom=1.5in]{geometry}
\usepackage{graphicx}
\usepackage{color}
\usepackage{amsmath, amssymb, amsthm,extarrows}
\definecolor{red}{rgb}{0.7,0,0}
\definecolor{grey}{RGB}{112,112,112}
\definecolor{blue}{RGB}{034,113,179}
\usepackage[colorlinks=true,citecolor=blue,linkcolor=red,urlcolor=grey,backref=page]{hyperref}

\newcommand{\koniec}{\begin{flushright}  $\Box $ \end{flushright}}
\newtheorem{theo}{Theorem}[section] 
\newtheorem{prop}[theo]{Proposition}  

\newtheorem{defi}[theo]{Definition}

\theoremstyle{remark}

\newcounter{mnotecount}[section]

\renewcommand{\themnotecount}{\thesection.\arabic{mnotecount}}

\newcommand{\mnote}[1]%{}%
{\protect{\stepcounter{mnotecount}}$^{\mbox{\footnotesize
$%\!\!\!\!\!\!\,
\bullet$\themnotecount}}$ \marginpar{%\color{red}%
\raggedright\tiny\em
$\!\!\!\!\!\!\,\bullet$\themnotecount: #1} }

\newcommand{\hook}{{\setlength{\unitlength}{11pt}   % adjust pt size here
                   \begin{picture}(.833,.8)
                   \put(.15,.08){\line(1,0){.35}}
                   \put(.5,.08){\line(0,1){.5}}
                   \end{picture}}}
\newcommand{\CP}{\mathbb{CP}}
\newcommand{\C}{\mathbb{C}}

\newcommand{\PP}{\mathbb{P}}
\newcommand{\RP}{\mathbb{RP}}
\newcommand{\tL}{\widetilde{L}}
\newcommand{\tP}{\widetilde{P}}
\newcommand{\R}{\mathbb{R}}

\def\p{\partial}
\def\be{\begin{equation}}

\def\ee{\end{equation}}
\def\ep{{\varepsilon}}
\def\bea{\begin{eqnarray}}
\def\eea{\end{eqnarray}}

\newcommand{\spp}{\mathbb{S}}

\numberwithin{equation}{section}
\begin{document} \date{29 August 2019}
%%%%%%%%%%%%%%%%%%%%%%%%%%%%%%%%%%%%%%%%
\title{ Einstein Metrics, Projective Structures and the $SU(\infty)$ 
Toda equation }

\author{Maciej Dunajski}
\address{Department of Applied Mathematics and Theoretical Physics\\ 
University of Cambridge\\ Wilberforce Road, Cambridge CB3 0WA, UK.}
\email{m.dunajski@damtp.cam.ac.uk}

\author{Alice Waterhouse}
\address{Department of Applied Mathematics and Theoretical Physics\\ 
University of Cambridge\\ Wilberforce Road, Cambridge CB3 0WA, UK.}
\email{aw592@cam.ac.uk}

\begin{abstract}
We establish an explicit correspondence between two--dimensional projective structures admitting a projective vector field, and 
a class of solutions to the $SU(\infty)$ Toda equation. We give several examples of new, explicit solutions of the
Toda equation, and construct their mini--twistor spaces. Finally we discuss the projective-to-Einstein correspondence, which gives
a neutral signature Einstein metric on a cotangent bundle  $T^*N$ of any projective structure $(N, [\nabla])$. We show that
there is a canonical Einstein of metric on an $\R^*$--bundle over $T^*N$, with a connection whose curvature is the pull--back
of the natural symplectic structure from $T^*N$. 
\end{abstract}
\maketitle
\begin{center} 
{\em Dedicated to Joseph Krasil'shchik on the occasion of his 70th birthday.}
\end{center}
\section{Introduction}
The aim of this paper is two--fold:
\begin{itemize}
\item[(A)] To associate a Lorentzian Einstein--Weyl structure in $(2+1)$--dimensions with any projective
structure on a surface which admits a one--parameter group of projective symmetries.
\item[(B)] To construct an explicit class of solutions of the $SU(\infty)$--Toda equation
\be
\label{md_toda}
U_{XX}+U_{YY}=\epsilon(e^U)_{ZZ}, \quad\mbox{where}\quad U=U(X, Y, Z), \quad
\mbox{and}\;\;\epsilon=\pm 1
\ee
with no continuous group of point symmetries.
\end{itemize}
We shall see that (B) is an explicit coordinate realisation of (A), but we 
chose to separate the two constructions
for the benefit of readers interested in integrable systems and solutions to (\ref{md_toda}) but not necessarily willing to study the relationships between the projective, conformal, and Weyl differential structures.

Equation (\ref{md_toda}) has originally arisen in  the context of complex general relativity \cite{FP, BF, Prz}, and then
in Einstein--Weyl \cite{ward_toda} and (in Riemannian context, with
$\epsilon=-1$) scalar--flat K\"ahler geometry \cite{LeBrun}. It belongs to a class
of dispersionless systems integrable by the twistor transform 
\cite{MW, MDbook, ADM}, 
the method of  hydrodynamic reduction \cite{F},  and  the Manakov--Santini approach \cite{MS}. 
The equation
is nevertheless not linearisable and most known explicit solutions admit Lie point or other symmetries (there are exceptions - see 
\cite{c_toda, CT,martina, Sheftel}, as well as \cite{K1, K2} where other general frameworks are discussed). 
The solutions
we find depend on two arbitrary functions of one variable, and arise from an essentially linear procedure, where no non--linear PDEs/ODEs have to be solved. An example of a solution in our class is given by an implicit relation
\be
\label{example_int}
4Y^2e^U(e^UX^2-Z^2)^3+(2e^{2U}X^4-3e^UX^2Z^2+Z^4+2Z^2)^2=0,
\ee
where the level sets of $U$ in $\R^3$ are real algebraic surfaces.

Now we move on to describe the construction (A), which is based on a combination of the Jones--Tod correspondence
\cite{JT}, a theorem of Tod \cite{Tod_note} which improved an earlier result of Przanowski \cite{Prz}, and two theorems from \cite{DM}.
In \cite{DM} it was shown that with any projective structure $[\nabla]$ on a surface $N$ one can associate
a neutral signature Einstein metric with non--zero scalar--curvature, and an anti--self--dual (ASD) Weyl tensor. If the projective
structure is represented by an affine connection $\nabla\in[\nabla]$ (see \S\ref{s_defini} for definitions), then the metric is isometric to the following metric on $T^*N$:
\be
\label{g}
g= dz_A \odot dx^A - (\Gamma^C_{AB}z_C-z_Az_B-P_{AB})dx^A \odot dx^B, 
\quad A, B, C=0, 1 , 
\ee
where $x^A$ are coordinates on $N$, $z_A$ are coordinates on the fibres of $T^*N$, $\Gamma^C_{AB}$ are connection components of $\nabla$, and $P_{AB}$ is the projective Schouten tensor of $\nabla$. The following Theorem has  been established 
in \cite{DM} 
\begin{theo}\cite{DM}
\label{th1int}
Let $(M, g)$ be an ASD Einstein manifold with scalar curvature 24 admitting a parallel ASD totally null distribution.  
\begin{enumerate}
\item
Then $(M, g)$ is conformally flat, or it is locally isometric to 
(\ref{g}) for
some torsion-free connection $\nabla$ on a surface $N$.
\item There is a one--to--one correspondence between projective vector fields of $(N, \nabla)$, and Killing vector fields of $(M, g)$ with
the metric $g$ given by (\ref{g}).
\end{enumerate}
\end{theo}
If the projective strucure $(N, [\nabla])$ admits a projective vector field, then the connection with
the Einstein--Weyl geometry now follows from Theorem \ref{th1int}, 
and the neutral signature
version of the Jones--Tod correspondence:
\begin{theo}\cite{JT}
\label{theo_tod1}
Let $(M, g)$ be a four--manifold with a neutral signature metric with ASD Weyl tensor, and a conformal Killing vector $K$. Let
\be 
\label{EWgen}
h=|K|^{-2}g-|K|^{-4}{\bf{K}}\odot{\bf{K}},\qquad \omega=\frac{2}{|K|^2}\star({\bf{K}}\wedge d{\bf{K}}),
\ee
where $|K|^2=g(K,K)$, ${\bf{K}}=g(K, \cdot)$ and $\star$ is the Hodge operator defined by $g$. Then $(h, \omega)$ is a solution of the Einstein--Weyl equations  on the space of orbits $W$ of $K$ in $M$. All Lorentzian Einstein--Weyl structures
arise from some anti-self-dual $(M, g, K)$.
\end{theo}
The final step to realising (B) is the occurrence of the $SU(\infty)$--Toda equation 
(\ref{md_toda}). This is a consequence of the following result
of Tod
\begin{theo}
\label{th3int}\cite{Tod_note}
Let $(h, \omega)$ be the Einstein--Weyl structure arising from Theorem \ref{theo_tod1}, under the additional assumption that
the ASD conformal structure $(M, g)$ is Einstein, and with non--zero Ricci scalar. 
\begin{enumerate}
\item The Einsten--Weyl structure admits a shear--free, twist--free geodesic
congruence.
\item There exists $h\in [h]$, and
coordinates $(X, Y, Z)$ on an open set in $W$ such that
(assuming the signature of $h$ is $(2, 1)$ and the congruence is time--like)
\be
\label{metric_toda}
h=e^U(dX^2+dY^2)-dZ^2, \quad \omega =2U_ZdZ
\ee
and the function $U=U(X, Y, Z)$ satisfies the $SU(\infty)$--Toda equation
(\ref{md_toda}) 
with $\epsilon=1$.
\end{enumerate}
\end{theo}
The whole construction can now be summarised in the following diagram
\begin{eqnarray}
\label{diagram}
\text{Projective structure with symmetry} &\overset{\text{Thm}\;\ref{th1int}}
{\longrightarrow}& \text{ASD Einstein with symmetry}\nonumber\\
\downarrow & & \downarrow \text{Thm}\;\ref{theo_tod1}\\
\text{Solution to}\;SU(\infty)\; \text{Toda} &\overset{\text{Thm}\;\ref{th3int}}\longleftarrow&\text{Einstein--Weyl}\nonumber
\end{eqnarray}
%\be
%\label{diagram}
%A \xrightarrow{\text{Theorem}} \text{ASD Einstein with symmetry}\xrightarrow{\t%ext{Theorem}} \text{Einstein--Weyl} \xrightarrow{\text{Theorem}} SU(\infty)\; \%text{Toda}
%\ee
The paper is organised as follows. In the next section we summarise the basic facts and relevant formulae underlying Theorems
\ref{th1int} and \ref{theo_tod1}. In proposition \ref{prop1} of \S \ref{general} we present the most general class of EW spaces arising from our construction, and in the remainder of the section we show how to associate  solutions of the $SU(\infty)$--Toda equation with this class. In \S\ref{neat} 
and \S\ref{neat2}
we give several examples corresponding to $SL(2, \R)$ and $SL(3, \R)$
invariant projective structures. In the latter case the four--manifold $(M, g)$ of Theorem \ref{th1int} is $SL(3)/GL(2)$, and the mini--twistor space
of the $SU(\infty)$--Toda equation can be constructed explicitly by quotienting the flag manifold $F_{12}(\C^3)$ by a $\C^*$ action. 
In Proposition \ref{ewsymprop} we give an explicit criterion, in terms of the representative metric $h\in [h]$ and the one form $\omega$ for a vector field the generate a symmetry of the Weyl structure.
In Proposition \ref{propHH} we show that the Einstein metric on $SL(3)/GL(2)$  is also pseudo--hyper--Hermitian, and its twistor space fibers holomorphically over $\CP^1$. In \S\ref{ODEs} we make contact with the Cartan approach to Einstein--Weyl geometry via
special 3rd order ODEs.
In \S \ref{KKsection} we shall prove (Theorem \ref{theokk}) that the 
$2n$--dimensional analogue of the Einstein metric (\ref{g}) canonically lifts
to an Einstein metric of signature $(n, n+1)$ on the $\R^*$ bundle 
${\mathcal Q}$ over $M$ with a connection whose curvature is the pullback of 
the symplectic form from $M$ to ${\mathcal Q}$. Some calculations
underlying the proof of Theorem \ref{theokk} are relegated to Appendix A.
In the Appendix B we shall present a solution to the elliptic $SU(\infty)$--Toda equation corresponding to an ALH gravitational instanton.
\subsection*{Acknowledgments} 
The work of MD has been partially supported by STFC consolidated grant no. ST/P000681/1. 
AW is grateful for support from the Sims Fund. We thank David Calderbank, Thomas Mettler and Jacek Tafel for useful discussions.
\section{Projective, Einstein, and Weyl geometries}
\label{s_defini}
Here we summarise basic facts about projective, Einstein--Weyl, 
and anti--self--dual geometries.
\subsection{Projective structures}
In this section we review projective differential geometry. In the applications to the $SU(\infty)$--Toda equation we 
shall focus on the surface case, where the dimension of the underlying manifold $N$ is two. In
\S\ref{KKsection} we shall consider the general case where $\mbox{dim}(N)=n$.
\begin{defi}
A projective structure on a surface $N$ is an equivalence class $[\nabla]$ of torsion-free affine connections on $TN$ which share the same unparametrised geodesics.
\end{defi}
Let $\nabla \in [\nabla]$ be a connection in the projective equivalence class
with connection symbols $\Gamma^{C}_{AB}$. Any other connection
in $[\nabla]$  can be obtained from $\nabla$ in terms of a one-form $\Upsilon$ as
\be
\label{proj_change}
\Gamma^{C}_{AB} \rightarrow \Gamma^{C}_{AB} + \delta^C_A\Upsilon_B + \delta^C_B\Upsilon_A, \quad A, B, C=1, 2, \dots, n. 
\ee
Let $R_{AB}$ be the (not necessarily symmetric) Ricci tensor of $\nabla$, and
let $P_{AB}=1/(n-1)R_{(AB)}+1/(n+1) R_{[AB]}$ be the projective Schouten tensor
(here we give the formula for the general $n$--dimensional projective structure).
The change of representative connection (\ref{proj_change}) induces
the following change
\be
\label{shouten_change}
P_{AB} \rightarrow P_{AB} + \Upsilon_A\Upsilon_B - \nabla_A\Upsilon_B 
\ee
to the Schouten tensor. A two--dimensional 
projective structure is called flat if it is locally diffeomorphic to the real projective plane
with unparametrised geodesics given by projective lines. This happens if and only if
the Cotton tensor $\nabla_{[A}P_{B]C}$ vanishes for any choice of the representative connection.

Let $(N, \nabla)$ be a manifold with an affine connection. A projective vector
field $k$ is a generator of a one--parameter group of transformations
mapping unparametrised geodesics of $\nabla$ to unparametrised geodesics.
At the infinitesimal level the projective condition is
\be 
\label{projvf}
\mathcal{L}_k\Gamma^C_{AB}=\delta^C_A\Upsilon_B+\delta^C_B\Upsilon_A,
\ee 
where the Lie derivate of the connection components is defined as in \cite{yano}.
In general no projective vector fields exist on $(N, \nabla)$. The possible Lie algebras of projective vector fields on a surface are $\mathfrak{sl}(3, \R)$, 
$\mathfrak{sl}(2, \R)$, $\mathfrak{a}_2$ (the two--dimensional affine Lie algebra) or $\R$. See \cite{Bryant} for further details.
\subsection{Einstein--Weyl structures}
\begin{defi}
A Weyl Structure $(W,\mathcal{D},[h])$ is a conformal equivalence class of metrics $[h]$ on a manifold $W$ along with a fixed torsion--free affine connection $\mathcal{D}$ which preserves any representative $h\in[h]$ up to conformal class. That is, for some one-form $\omega$,
\[
\mathcal{D}h=\omega\otimes h.
\]
\end{defi}
A pair $(h,\omega)$ uniquely defines the connection and hence the Weyl structure, but there is an equivalence class of such pairs which define the same Weyl structure. These are related by transformations
\be
\label{weyl_tr}
h\rightarrow \rho^2h,\quad\omega\rightarrow\omega+2d\mathrm{ln}(\rho),
\ee
where $\rho$ is a smooth, non-zero function on $W$. 

If additionally the symmetric part of the Ricci tensor of $\mathcal{D}$ is a scalar multiple of $h$, then $W$ is said to carry an Einstein-Weyl structure. Physically, the Einstein--Weyl condition in Lorentzian signature corresponds to the statement that null geodesics of the conformal structure $[h]$ are also geodesics of the connection $\mathcal{D}$.
This condition is invariant under (\ref{weyl_tr}). In three dimenstions, the Einstein--Weyl equations
give a set of five non--linear PDEs on the pair $(h, \omega)$. These equations
are integrable by the twistor transform of Hitchin \cite{hitchin}, which 
(by Theorem \ref{theo_tod1}) can be regarded as a reduction of Penrose's twistor transform \cite{penrose} for ASD conformal structures.
A trivial Einstein--Weyl structure is one whose one--form $\omega$ is closed, so that it is locally exact and thus may be set to zero by a change of scale (\ref{weyl_tr}). Then $\mathcal{D}$ is the Levi--Civita connection of some representative $h\in[h]$, and this representative is Einstein.
\subsection{Anti--self--dual Einstein metrics}
Let $M$ be an oriented four--dimensional manifold with a metric
$g$ of signature $(2, 2)$. The Hodge $\ast$ operator on the space of two forms is an involution, and induces a decomposition \cite{AHS}
\be 
\label{split_ad}
\Lambda^{2}(T^*M) = \Lambda_{-}^{2}(T^*M) \oplus \Lambda_{+}^{2}(T^*M)
\ee
of two-forms
into anti-self-dual (ASD)
and self-dual (SD)  components, which
only depends on the conformal class of $g$. 
The Riemann tensor of $g$ 
can be thought of
as a map $\mathcal{R}: \Lambda^{2}(T^*M) \rightarrow \Lambda^{2}(T^*M)$
which admits a decomposition   under (\ref{split_ad}):
\be \label{decomp}
{\mathcal R}=
\left(
\mbox{
\begin{tabular}{c|c}
&\\
$C_+-2\Lambda$&$\phi$\\ &\\
\cline{1-2}&\\
$\phi$ & $C_--2\Lambda$\\&\\
\end{tabular}
} \right) ,
\ee
where $C_{\pm}$ are the SD and ASD parts 
of the Weyl tensor, $\phi$ is  the
trace-free Ricci curvature, and $-24\Lambda$ is the scalar curvature which acts
by scalar multiplication. 
The metric $g$ is ASD and Einstein if 
$C_{+}=0$ and $\phi=0$. In this case
the Riemann tensor is also anti-self-dual.
\vskip3pt
Locally there exist real rank-two vector bundles $\spp, \spp'$  called spin-bundles over $M$, equipped with parallel symplectic structures
$\ep, \ep'$ such that
\be
\label{can_bun_iso}
T M\cong {\spp}\otimes {\spp'}
\ee
is a  canonical bundle isomorphism, and
\[
g(v_1\otimes w_1,v_2\otimes w_2)
=\varepsilon(v_1,v_2)\varepsilon'(w_1, w_2)
\]
for $v_1, v_2\in \Gamma(\spp)$ and $w_1, w_2\in \Gamma(\spp')$.
A vector $V\in \Gamma(TM)$ is called null if $g(V, V)=0$. Any null vector is of the form
$V=\lambda \otimes \pi$ where $\lambda$ and $\pi$ are sections of
$\spp$ and $\spp'$ respectively.
An $\alpha$--plane (respectively a $\beta$--plane) 
is a two--dimensional plane in $T_pM$
spanned by null vectors of the above form with $\pi$ (respectively
$\lambda$) fixed, and
an $\alpha$--surface ($\beta$--surface) is a two--dimensional surface in $M$ 
such 
that its tangent plane at every point is an $\alpha$--plane ($\beta$--plane). 
Penrose's Nonlinear Graviton Theorem \cite{penrose} states that   a 
maximal three dimensional family of $\alpha$--surfaces exists 
in $M$ iff $C_+=0$.
\subsubsection{ASD Einstein metrics from projective structures}
A general ASD metric depends, in the real--analytic category, on six arbitrary functions of three variables. Theorem \ref{th1int} gives an explicit subclass of such metrics which are additionaly Einstein and carry a so--called parallel ASD totally null distribution. These depend on two arbitrary functions of two variables. Any projective structure $(N, [\nabla])$ gives rise to such an ASD Einstein metric. The explicit form (\ref{g}) is the pull-back of the metric on $M$ along the diffeomorphism $\varphi:T^*N\rightarrow M$ specified by the choice of connection $\nabla\in [\nabla]$. In \cite{DGW} it is shown how to extend this metric to a $c$--projective compactification.

There is an additional structure on four--manifolds described in Theorem \ref{th1int}: a para--Hermitian structure. The symplectic form $\Omega$ of this para--Hermitian structure pulls back to
\be
\label{Omega}
\varphi^{*}{\Omega} = dz_A\wedge dx^A + P_{AB}dx^A \wedge dx^B,
\ee
where $(x^A,z_A)$ are canonical local coordinates on the cotangent bundle.
The pair $(g, \Omega)$ is projectively invariant under the changes 
(\ref{projvf}) if $z_A\rightarrow z_A+\Upsilon_A$.

If $k$ is a projective vector field on $(N, \nabla)$, then the corresponding Killing vector field on $(M, g)$ is symplectic, and 
is given in local coordinates by
\be 
\label{symlift}
K=k-z_A\frac{\p k^B}{\p x^A}\frac{\p}{\p z_B} + \Upsilon_A\frac{\p}{\p z_A}.
\ee
\subsubsection{ASD $\beta$--foliation}
It follows from the general construction of Calderbank \cite{Cal1}
and West \cite{West}
that any ASD conformal structure arising from Theorem \ref{th1int} 
carries a foliation by $\beta$--surfaces defined by an ASD two--form
$\Sigma_{ab}=\iota_A\iota_B\epsilon_{A'B'}$, and  such that
the spinor $\iota_A$ satisfies
\be
\label{dm3}
\nabla_{A'(A}\iota_{B)}=\mathcal{A}_{A'(A}\iota_{B)}
\ee
where $d\mathcal{A}$ is an ASD Maxwell field.

We shall call such foliations ASD $\beta$--surface foliations. In our coordinates $\Sigma=dx^0\wedge dx^1$ and ${D}=\mathrm{span}\{\p/\p z_0,\p/\p z_1\}$. We find that
\be
\label{beta_eq}
\nabla\Sigma=6\mathcal{A}\otimes \Sigma,
\ee
where $d\mathcal{A}=\Omega$, and $\Omega$ is the symplectic form on $M$, whose anti-self-duality implies (\ref{dm3}) for a rescaling of $\mathcal{A}$. In \S\ref{model} we will consider the model case where $M=SL(3)/GL(2)$, for which we can find the Ward transform of this ASD Maxwell two-form to the twistor space $F_{12}(\C^3)$
explicitly, and we will find that there is a second ASD $\beta$--surface foliation with a different two--form.

\section{From projective structures to $SU(\infty)$ Toda fields}
\label{general}
Recall (see, e.g. \cite{BDE}) that a projective structure on a surface can be locally specified by a single 2nd order ODE: taking coordinates $(x,y)$ on the surface we find that geodesics on which $\dot{x}\neq 0$ can be written as unparametrised curves $y(x)$ such that
\be
\label{odealice}
y^{\prime \prime} + a_0(x,y)+3a_1(x,y)y^{\prime}+3a_2(x,y)(y^{\prime})^2 + a_3(x,y)(y^\prime)^3=0,
\ee
where the coefficients $\{a_i\}$ are given by the projectively invariant formulae
\[
a_0=\Gamma^1_{00},\quad
3a_1=-\Gamma^0_{00}+2\Gamma^1_{01},\quad
3a_2=-2\Gamma^0_{01}+\Gamma^1_{11},\quad
a_3=-\Gamma^0_{11}.
\]
Consider the most general Einstein--Weyl structure arising from the combination
of Theorem \ref{th1int} and Theorem \ref{theo_tod1}. Because of the correspondence 
(Theorem \ref{th1int}, part 2.)
between symmetries of $(M, g)$ and symmetries of the projective surface $(N, [\nabla])$, the construction must begin with the general projective surface with at least one symmetry. In this case, the unparametrised geodesics can generically be written uniquely as integral curves of
the ODE
\be
\label{1symODE}
y^{\prime \prime}=A(y)(y^{\prime})^3+B(y)(y^{\prime})^2+1.
\ee
This has been established in \cite{FLL}, and the argument justifying this normal form is as follows.
We can locally choose the coordinates so that the projective symmetry is generated by the translation
$k=\p/\p x$. The normalising coordinates for $\frac{\p}{\p x}$ are unique up to a change
$(x, y) \rightarrow (x + \phi(y), \psi(y)).$
As long as $a_0$ is not everywhere zero (the generic case), we can choose $a_0=-1$ using the freedom
$(x, y) \rightarrow (x, \psi(y))$,
by taking $\psi'=-1/a_0$. This removes the $\psi$ freedom up to an additive
constant. Then we can choose $a_1=0$ using a change
$(x, y) -> (x + \phi(y), y)$,
where $\phi'=-a_1$. This preserves $a_0=-1$, and removes the $\phi$ freedom up to an
additive constant.

The  projective structure resulting from \ref{1symODE} is flat iff both
$A$ and $B$ are constant.
By trial and error, we chose a representative connection for (\ref{1symODE}) such that (\ref{g}) had the simplest possible form. The choice of connection we took was
\[
\Gamma^{0}_{11}=A(y),\quad \Gamma^{1}_{00}=-1, \quad \Gamma^{1}_{11}=-B(y)
\]
with all other components vanishing. Note that this choice of connection has a symmetric Ricci tensor, so the Schouten tensor is also symmetric and the symplectic form (\ref{Omega}) pulls back to just $dz_A\wedge dx^A$. Thus we can write the Maxwell potential $\mathcal{A}$ which is such that $d\mathcal{A}=\Omega$ as $\mathcal{A}=z_Adx^A$. Writing $x^A=(x,y)$, $z_A=(p,q)$, the resulting metric (\ref{g}) is
\be
\label{einstein_1}
g=(B(y) + p^2 +q)dx^2+2(pq+A(y))dxdy + (-A(y)p+B(y)q+q^2)dy^2 + dxdp +dydq.
\ee
Factoring by $K=\frac{\p}{\p x}$, and following the algorithm of Theorem
\ref{theo_tod1} gives the following
\begin{prop}
\label{prop1}
The most general  Einstein--Weyl structure arising
from the procedure (\ref{diagram}) is locally equivalent to
\begin{eqnarray}
\label{ew_final}
h&=&\frac{1}{V}\big((Bq -Ap+ q^2)dy+dq\big)dy
-\Big((pq+A)dy+\frac{1}{2}dp\Big)^2, \label{genh} \\
\omega&=&V(4dq+2 pdp), \qquad V=
({B}+ p^2+q)^{-1},\nonumber
\end{eqnarray}
where $(p, q, y)$ are local coordinates on $W$, and $A, B$ are arbitrary functions of $y$.
\end{prop}
\subsection{Solution to the $SU(\infty)$--Toda equation}
\label{steps_sec}
The procedure for extracting the corresponding solution to the $SU(\infty)$--Toda equation is given in \cite{Tod_note} 
(see also \cite{LeBrun} and \cite{DTkahler}). It involves finding the coordinates $(X,Y,Z)$ that put the metric (\ref{genh}) in the 
form (\ref{metric_toda}). Given an ASD Einstein metric $(M, g)$ with a
Killing vector $K$
\begin{enumerate}
\item The conformal factor $c:M\rightarrow \R^+$  given by
\[
c={|d{\bf K}+*_g d{\bf K}|_{g}}^{-1/2}
\]
has a property that
the rescaled self--dual derivative of $K$
\[
\Theta\equiv c^3\Big(\frac{1}{2}(d{\bf K}+*_g d{\bf K})\Big)
\]
is parallel with respect to $\hat{g}=c^2 g$.
The metric $\hat{g}$ is K\"ahler with  self--dual 
K\"ahler form $\Theta$, and admits a Killing vector $K$, as
${\mathcal L}_K(c)=0$.
\item
Define a function $Z:M\rightarrow \R$ to be the moment map:
\be
\label{ztilde}
dZ=K\hook \Theta.
\ee
It is well defined, as the K\"ahler form is Lie--derived along $K$.
\item
Construct the Einstein--Weyl structure of Theorem \ref{theo_tod1}
by factoring $(M, \hat{g})$ by $K$. Restrict the metric $h$
to a surface $Z=Z_0=\mbox{const}$, and construct isothermal coordinates 
$(X, Y)$ on this surface:
\[
\gamma\equiv h|_{Z=Z_0}=e^{U}(dX^2+dY^2), \quad U=U(X, Y, Z_0).
\]
To implement this step chose an orthonormal basis of one--forms
such that $\gamma= {e_1}^2+{e_2}^2$. Now $(X, Y)$ are solutions to the linear
system of 1st order PDEs
\[
(e_1+ie_2)\wedge (dX+idY)=0.
\]
\item Extend the coordinates $(X, Y)$ from the surface $Z=Z_0$ to $W$. This may
involve a $Z$--dependent affine transformation of $(X, Y)$.
\end{enumerate}
Implementing the steps 1--4 on MAPLE we find that if $A=0$, and $B=B(y)$ 
is arbitrary, then the $SU(\infty)$--Toda solution is given implicitly by
\begin{eqnarray}
\label{toda_implicit1}
 X&=&-\frac{8\mathrm{e}^{-2\int{B(y)dy}}Z^3p}{(Z^2p^2+4)^2},\quad
Y=\int{\mathrm{e}^{-2\int{B(y)dy}}dy}+\frac{\mathrm{e}^{-2\int{B(y)dy}}(-2Z^4p^2+8Z^2)}{(Z^2p^2+4)^2}.\nonumber\\
U&=&\mathrm{ln}\bigg(\frac{(Z^2p^2+4)^3}{64Z^2}\bigg)+4\int{B(y)dy}.
\end{eqnarray}
We can check that this is indeed a solution using the fact that the $SU(\infty)$--Toda equation is equivalent to 
$d\star_h dU=0$. We have also checked by performing a coordinate transformation of (\ref{md_toda}) to the coordinates $(y,p,Z)$.
\subsubsection{Example 1.}
Consider the flat projective structure with $A=B=0$, in which case
the coordinate $p$ can be eliminated between
\[
e^U=\bigg(\frac{(Z^2p^2+4)^3}{64Z^2}\bigg), \quad
X=-\frac{8Z^3p}{(Z^2p^2+4)^2}
\]
by taking a resultant. This yields
%\[
%\mathrm{e}^{4U}X^6+3\mathrm{e}^{3U}X^4Z^2+3\mathrm{e}^{2U}X^2Z^4+\mathrm{e}^UZ^%6+Z^4=0.
%\]
\[
e^U(e^UX^2-Z^2)^3+Z^4=0.
\]
\subsubsection{Example 2.}
To simplify the form of (\ref{toda_implicit1}) set
\[
G=\int\exp{\Big(-2\int B(y)dy\Big)}, \quad T=\frac{2Z^2}{Z^2p^2+4}
\]
Then (\ref{toda_implicit1}) becomes
\be
\label{toda_impl2}
e^U=\frac{Z^4}{8T^3 (G')^2}, \quad Y=G+G'T\Big(\frac{4T}{Z^2}-1\Big), \quad
X^2=\frac{4T^4(G')^2}{Z^2}\Big(\frac{2}{T}-\frac{4}{Z^2} \Big).
\ee
Eliminating $(T, y)$ between these three equations gives one relation between $(X, Y, Z)$ and  $U$ - this is our implicit solution.
The elimination can be carried over explicitly if $G=y^k$ for any integer $k$, or if $G=\exp{y}$. In the later case the solution is given by
(\ref{example_int}).
\subsection{Two monopoles}
The Einstein--Weyl structures  (\ref{genh}) we have constructed in Proposition \ref{prop1}
are special, as they belong to the $SU(\infty)$--Toda class, and so
(as shown by Tod \cite{Tod_toda}) admit a non--null geodesic congruence which has vanishing shear and twist.
The general solution to the $SU(\infty)$--Toda equation depends (in the real analytic category) on 
two arbitrary functions of two variables, but the solutions of the form (\ref{genh}) depend on two functions of one variable. The additional constraints on the solutions can be traced back to the four dimensional ASD conformal structures
which give rise (by the Jones--Tod construction) to (\ref{genh}). In addition to their being ASD and Einstein
they are  characterised \cite{DM} by a $\beta$--distribution which is parallel with respect to the Levi--Civita connection and ASD in the sense of Calderbank \cite{Cal1}.
The corresponding $\beta$--surfaces do not generically intersect with a given $\alpha$--surface, however if they do intersect then they will intersect in curves (null geodesics) which descend to the Einstein--Weyl structures, and give
rise to another (in addition to the Tod shear--free, twist--free) geodesic congruence.
In what follows we shall point out how some of this structure arises from a couple of solutions to the Abelian monopole
equation on EW  backgrounds.

Jones and Tod \cite{JT} show that there is a correspondence between conformally ASD four-metrics over an Einstein-Weyl structure $(W,\mathcal{D},[h])$, and solutions to the monopole equation on $W$. Hence the symmetry reduction of Theorem
\ref{theo_tod1} gives us a solution to the monopole equation. In fact, since we have an ASD Maxwell field on $M$, the reduction gives us a second monopole. In this subsection we compute these explicitly.
Given an EW structure $(h, \omega)$ in 2+1 dimension, the Abelian monopole consits
of a pair $(V, \alpha)$, where $V$ is a function, and $\alpha$ is a one--form subject to the equation
\[
dV+\frac{1}{2}\omega V=\star_h d\alpha.
\]
The inverse Jones--Tod correspondence \cite{JT} associates a neutral signature 
ASD conformal structure
\[
g=Vh-V^{-1}(dx +\alpha)^2
\]
with an isometry $K=\p/\p x$ with any solution of the monopole equation.

The conformal gauge in the EW geometry of Proposition \ref{prop1} 
is chosen so that 
\[
\alpha=V( pq+A)dy+\frac{V}{2}dp.
\]
Let us call this solution the Einstein monopole, as the resulting conformal class contains an Einstein metric
(\ref{einstein_1}). The second solution $(V_M, \alpha_M)$ (which we shall call the Maxwell monopole)
arises
as a symmetry reduction of the ASD Maxwell potential
\[
{\mathcal A}=pdx+qdy=-V_M K+\alpha_V,
\]
where $K=K_{\mu}dx^{\mu}$ is the Killing one--form, and
we find
\[
V_M=-pV, \quad \alpha_V=qdy-p\alpha.
\]
\section{An example from the submaximally symmetric projective surface}
\label{neat}
The submaximally symmetric projective surface is the punctured plane $N=\mathbb{R}^2\backslash{0}$, with the symmetry group $SL(2)$ acting via its fundamental representation. Here we have a one parameter family of projective structures falling into three distinct equivalence classes, with geodesics described by the differential equation
\be
\label{neatODE}
y^{\prime \prime}=-\mu(y-xy^\prime)^3,
\ee
where $(x,y)$ are coordinates on $\mathbb{R} ^2$ and $\mu$ is a constant parameter. The equivalence class that a given projective structure falls into depends on the value of $\mu$: those with $\mu>0$ form one of the classes, those with $\mu<0$ form another, and those with $\mu=0$ form the third. Further details can be found in \cite{Bryant}. For simplicity, we choose $\mu=1$.

Choosing a representative connection from the projective class defined by (\ref{neatODE}), we obtain from (\ref{g}) an
Einstein metric
\be
\begin{split}
g=( p^2- xy^2p-y^3q + 4 y^2 )dx^2 + 2(pq +  x^2yp +  xy^2q - 4 xy)dxdy \\
+ (q^2 - x^3p -  x^2yq + 4 x^2)dy^2 + dxdp + dydq
\end{split}
\ee
on $M$, again with $z_0=:p,\,z_1=:q$, having Killing vectors
\[
K_1=x\frac{\p}{\p x} - p\frac{\p}{\p p} - y\frac{\p}{\p y} + q\frac{\p}{\p q},\quad
K_2=x\frac{\p}{\p y} - q\frac{\p}{\p p}, \quad
K_3=y\frac{\p}{\p x} - p\frac{\p}{\p q}.
\]
These are lifts of the projective vector fields corresponding to the $\mathfrak{sl}(2)$ elements
\[
T_1=\begin{pmatrix}\epsilon & 0\\
0 & -\epsilon
\end{pmatrix}
\quad
T_2 = \begin{pmatrix}0 & 0\\
\epsilon & 0
\end{pmatrix}
\quad
T_3 = \begin{pmatrix} 0 & \epsilon\\
0 & 0
\end{pmatrix}.
\]

Factoring by $K_3$ and choosing coordinates
\[
u=\frac{p^2}{ y^2}, \quad v=2\ln( y^2), \quad w=xp+yq,
\]
we obtain an Einstein-Weyl structure
\begin{eqnarray}
\label{ew_neat}
h&=&-du^2-2dudw-w(w^2+u-5w+4)dv^2+2(u-w+4)dvdw,\\
\omega&=&\frac{1}{u-w+4}du-\frac{3w}{u-w+4}dv-\frac{4}{u-w+4}dw.\nonumber
\end{eqnarray}
The solution to the $SU(\infty)$--Toda equation (\ref{md_toda}) which determines the Einstein-Weyl structure (\ref{ew_neat}) is described by an algebraic curve $f(\mathrm{e}^U,X,Y,Z)=0$ of degree six in $\mathrm{e}^U$ and degree twelve in the other coordinates. 
This solution has been found following th Steps 1-4 in \S\ref{steps_sec}, 
and is given by
\be
\begin{split}
64\mathrm{e}^{6U}X^6(X+Y)^3(X-Y)^3
-92\mathrm{e}^{5U}X^4Z^2(X+Y)^3(X-Y)^3 \\
+48\mathrm{e}^{4U}X^2Z^2(5X^6Z^2-14X^4Y^2Z^2+13X^2Y^2Z^2-4Y^4Z^2+9X^4+27X^2) \\
+8\mathrm{e}^{3U} Z^4(-20X^6Z^2+48X^4Y^2Z^2-36X^2Y^4Z^2+8Y^6Z^2-81X^4-243X^2Y^2) \\
+3\mathrm{e}^{2U}Z^4
(20X^4Z^4-36X^2Y^2Z^4+16Y^4Z^4+108X^2Z^2+216Y^2Z^2+243) \\
+6\mathrm{e}^UZ^8(-2X^2Z^2+2Y^2Z^2-9) +Z^{12}\\ =0.
\end{split}
\nonumber
\ee
Note that the  formulae (\ref{ew_neat}) are independent of the coordinate $v$, and therefore have a symmetry. This was unexpected because there is no other symmetry of $(M,g)$ that commutes with $K_3$. However, it is possible for symmetries to appear in the Einstein-Weyl structure without a corresponding symmetry of the ASD conformal structure. This can be seen from the general formula (\ref{EWgen}); the function $V$ may depend on the coordinate $v$ so that $g$ depends on $v$ even though $h$ does not. For example, the Gibbons-Hawking metrics \cite{GH} give a trivial Einstein--Weyl structure with the maximal symmetry group, but the four-metric is in general not so symmetric. Our discovery of this unexpected symmetry motivated a more concrete description of a symmetry of a Weyl structure.

\begin{defi}
An infinitesimal symmetry of a Weyl structure $(W,\mathcal{D},[h])$ is a vector field $\mathcal{K}$ which is both an affine vector field with respect to the connection\footnote{Recall that an affine vector field of a connection $\mathcal{D}$ is one which preserves its components, i.e. $\mathcal{L}_\mathcal{K}\Gamma^i_{jk}=0$.} $\mathcal{D}$ and a conformal Killing vector with respect to the conformal structure $[h]$.
\end{defi}

\begin{prop}
\label{ewsymprop}
Given an infinitesimal symmetry $\mathcal{K}$ of a Weyl structure $(W,\mathcal{D},[h])$ in dimension $N$, and a representative $h\in[h]$ such that $\mathcal{D}h=\omega\otimes h$, there exists a smooth function $f:W\rightarrow \R$ such that
\be
\label{EWsym}
\mathcal{L}_\mathcal{K}h=fh,\qquad\mathcal{L}_\mathcal{K}\omega=\frac{1}{N}d[\mathcal{K}\hook d(\mathrm{ln}(\mathrm{det}(h)))].
\ee
\end{prop}
\noindent\textbf{Proof.} The first equation follows immediately from the fact that $\mathcal{K}$ is a conformal Killing vector of $h$. It remains to evaluate the Lie derivative of the one--form $\omega$ along the flow of $\mathcal{K}$ given that $\mathcal{L}_\mathcal{K}h=fh$ and $\mathcal{L}_\mathcal{K}\Gamma^i_{jk}=0$, where $\Gamma^i_{jk}$ are the components of the connection $\mathcal{D}$. We do this by considering the Lie derivative of $\mathcal{D}h$:
\begin{align*}
\mathcal{L}_\mathcal{K}(\mathcal{D}_ih_{jk}) &= \mathcal{L}_\mathcal{K}(\p_ih_{jk})-\mathcal{L}_\mathcal{K}(\Gamma^l_{ji}h_{lk}+\Gamma^l_{ki}h_{jl}) \\
&= \mathcal{L}_\mathcal{K}(\p_ih_{jk}) - f(\Gamma^l_{ji}h_{lk}+\Gamma^l_{ki}h_{jl}).
\end{align*}
Now
\begin{align*}
\mathcal{L}_\mathcal{K}(\p_ih_{jk}) &= \mathcal{K}^l\p_l\p_ih_{jk}+(\p_i\mathcal{K}^l)\p_lh_{jk}+(\p_j\mathcal{K}^l)\p_ih_{lk} + (\p_k\mathcal{K}^l)\p_ih_{jl} \\
&= \p_i[\mathcal{K}^l\p_lh_{jk}+(\p_j\mathcal{K}^l)h_{lk}+(\p_k\mathcal{K}^l)h_{jl}] - (\p_i\p_j\mathcal{K}^l)h_{lk} - (\p_i\p_k\mathcal{K}^l)h_{jl}.
\end{align*}
The term with square brackets is just
\[
\p_i(\mathcal{L}_\mathcal{K}h_{jk})=\p_i(fh_{jk})=f\p_ih_{jk}+\p_ifh_{jk},
\]
so we have
\[
\mathcal{L}_\mathcal{K}(\mathcal{D}_ih_{jk})=f\mathcal{D}_ih_{jk}+\p_ifh_{jk}- (\p_i\p_j\mathcal{K}^l)h_{lk} - (\p_i\p_k\mathcal{K}^l)h_{jl}.
\]
Setting this equal to $\mathcal{L}_\mathcal{K}(\omega_ih_{jk})=(\mathcal{L}_\mathcal{K}\omega_i)h_{jk}+f\omega_ih_{jk}$ and cancelling $f\omega_ih_{jk}$ with $f\mathcal{D}_ih_{jk}$, we find
\begin{eqnarray}
(\mathcal{L}_\mathcal{K}\omega_i)g_{jk} &=& \p_ifh_{jk} - (\p_i\p_j\mathcal{K}^l)h_{lk} - (\p_i\p_k\mathcal{K}^l)h_{jl}\nonumber \\
\label{liederivom}
\implies\ \mathcal{L}_\mathcal{K}\omega_i &=& \p_if - \frac{2}{N}\p_i\p_j\mathcal{K}^j.
\end{eqnarray}
Finally, we note that
\[
\p_i\p_j\mathcal{K}^j=\frac{N}{2}\p_if-\frac{1}{2}\p_i[\mathcal{K}\hook d(\mathrm{ln}(\mathrm{det}(h))].
\]
This follows from tracing the expression $\mathcal{L}_\mathcal{K}h_{ij}=fh_{ij}$:
\be
\begin{gathered}
\nonumber
\mathcal{L}_\mathcal{K}h_{ij} = \mathcal{K}^k\p_kh_{ij} + (\p_i\mathcal{K}^k)h_{kj} + (\p_j\mathcal{K}^k)h_{ik} = fh_{ij} \\
\implies \quad \mathcal{K}^kh^{ij}\p_kh_{ij} + 2\p_k\mathcal{K}^k = Nf \\
\implies \quad 2\p_i\p_k\mathcal{K}^k = N\p_if -  \p_i(\mathcal{K}^kh^{jl}\p_kh_{jl})
\end{gathered}
\ee
and recalling that $h^{jl}\p_kh_{jl}=\p_k\mathrm{ln}(\mathrm{det}(h))$.
Substituting into (\ref{liederivom}) then yields the result.
\begin{flushright}
$\square$
\par\end{flushright}

We can easily verify the invariance of (\ref{EWsym}) under Weyl transformations. Let $(\hat{h},\hat{\omega})$ be a new metric and one--form related to the old ones by (\ref{weyl_tr}). Then
\[
\mathcal{L}_\mathcal{K}\hat{\omega}=\mathcal{L}_\mathcal{K}\omega + 2\mathcal{K}\hook d\mathrm{ln}(\rho)
\]
from (\ref{weyl_tr}), and from (\ref{EWsym}) we have
\begin{align*}
\mathcal{L}_\mathcal{K}\hat{\omega} & = \frac{1}{N}d[\mathcal{K}\hook d(\mathrm{ln}(\rho^{2N}\mathrm{det}(h)))] \\
& = \frac{1}{N}d[\mathcal{K}\hook d(\mathrm{ln}(\mathrm{det}(h)))] + \frac{2N}{N}\mathcal{K}\hook d\mathrm{ln}(\rho) \\
& = \mathcal{L}_\mathcal{K}\omega + 2\mathcal{K}\hook d\mathrm{ln}(\rho),
\end{align*}
as above. Note that the function $f$ in (\ref{EWsym}) will change according to
\[
\hat{f}=f+2\mathcal{K}\hook d\mathrm{ln}\rho.
\]

In the case of the Weyl structure (\ref{ew_neat}), the infinitesimal symmetry is
\[
\mathcal{K}=\frac{\p}{\p v}.
\]
Since we have chosen a scale such that $\mathcal{K}$ is in fact a Killing vector of $h$, we have that $\mathcal{K}\hook d(\mathrm{ln}(\mathrm{det}(h))=0$, so the one--form $\omega$ is also preserved by $\mathcal{K}$. This is consistent with the fact that it has no explicit $v$--dependence.
\section{The model $SL(3)/GL(2)$ and its reductions}
\label{model}
In the following section we discuss the four--manifold $(M,g)$ obtained from the maximally symmetric flat 
projective surface $N=\RP^2$. In this case, $g$ is the indefinite analogue of the Fubini--Study metric, and is not only bi--Lagrangian but also para--K\"ahler, since the symplectic form $\Omega$ is parallel with respect to the Levi--Civita connection of $g$. Chosing a representative connection with
$\Gamma_{AB}^C=0$ gives $g$ as
\be
\label{special_ein}
g=dz_A\odot dx^A+ z_Az_B dx^A\odot dx^B.
\ee
We begin by discussing the  conformal structure
of (\ref{special_ein}), both explicitly and in terms of twistor lines. We then note some structure which is unique to the model case: hyper-hermiticity and a second foliation by $\beta$--surfaces which is ASD in the sense of Calderbank. 

Finally, we present a classification of the Einstein--Weyl structures which can be obtained by Jones--Tod factorisation of $SL(3)/GL(2)$ and exhibit an explicit example of such a factorisation from the twistor perspective, reconstructing the conformal structure on $W$ from minitwistor curves.
\subsection{Conformal Structure}
\label{model_conf}
Let $M\subset \PP^2\times {\PP^2}^*$ be set of non--incident pairs 
$(P, L)$, where $P\in \PP^2$, and $L\subset \PP^2$ is a line. 
\begin{prop}
\label{prop_cone}
Two pairs $(P, L)$ and $(\tP, \tL)$ are null--separated
with respect to the conformal structure (\ref{special_ein})
if there exists
a line which contains three points $(P, \tP, L\cap \tL)$. 
\end{prop}
{\bf Proof.}
The  null condition of Proposition \ref{prop_cone}
defines a co--dimension one cone in $TN$: 
generically there is no line through three given points.

To find an analytic expression for the resulting conformal structure 
consider two pairs  $(P, L)$ and $(\tP, \tL)$ 
of non--incident points and lines. Let $L+t\tL$ be a pencil of lines. 
There exists $t$ such that 
\be
\label{dm1}
P\cdot (L+t\tL)=0,  \quad
\tP\cdot (L+t\tL)=0.
\ee 
Eliminating $t$ from 
(\ref{dm1}) gives
\[
(P\cdot L)(\tP\cdot \tL)-(\tP\cdot L)(P\cdot \tL)=0.
\]
Setting $\tP=P+dP, \tL=L+dL$ yields a metric 
$g$  representing the conformal structure
\[
g=\frac{dP\cdot dL}{P\cdot L}-\frac{1}{(P\cdot L)^2}(L\cdot dP)(P\cdot dL).
\]
We can use the normalisation $P\cdot L=1$, so that $P\cdot dL=-L\cdot dP$,
and
\be
\label{dm_metric}
g={dP\cdot dL}+(L\cdot dP)^2.
\ee
We take affine coordinates 
\be
\label{DM_parameter}
P=[x^A, 1],\quad L=[z_A, 1-x^Az_A]
\ee 
with a normalisation $P\cdot L=1$ to recover the metric (\ref{g}) which now takes the form (\ref{special_ein}).
\koniec
\subsection{Twistor space}
\label{twist_SSS}
Let $F_{12}(\C^3)\in \PP^2\times {\PP^2}^*$ be set of incident pairs 
$(p, l)$, so that $p\cdot l=0$. This is the twistor space of $(M, g)$.
A $\PP^1$ embedding corresponding to a point $(P, L)$
consists of all lines $l$ thorough $P$, and all points
$p=l\cap L$:
\be
\label{dm2}
P\cdot l=0, \quad p\cdot L=0, \quad p\cdot l=0.
\ee
Let $(P, L)$ and $(\tP, \tL)$ be null separated. The corresponding lines in $F_{12}$
intersect at a point $(p, l)$ given by
\[
p=L\wedge \tL, \quad l=P\wedge \tP,
\]
where $[L\wedge\tL]^{\alpha}=\epsilon^{\alpha\beta\gamma}L_{\alpha}\tL_{\beta}$ etc.
The incidence condition $p\cdot l=0$ now gives the conformal structure
(\ref{dm1}). The contact structure on $F_{12}$ is $1/2(l\cdot dp-p\cdot dl)=p\cdot dl$.

We shall now give an explicit parametrisation of twistor lines, and show how 
the metric (\ref{dm_metric}) arises from the Penrose condition 
\cite{penrose, ward}.
Let $P\in \PP^2$. The corresponding $l\in {\PP^2}^*$ is
\[
l=P\wedge \pi, \quad \mbox{where}\quad  \pi\sim a\pi+b P,
\]
where $a\in \R^*, b\in \R$. Thus $\pi$ parametrises a projective line $\PP^1$,
and by making a choice of $b$ we can take
$
\pi=[\pi^0, \pi^1, 0], 
$ where $\pi^A=[\pi^0, \pi^1]\in \PP^1$. The constraint $P\cdot l=0$ now holds.
To satisfy the remaining constraints in (\ref{dm2}) we take
\[
p=L\wedge l=(L\cdot\pi)P-(L\cdot P)\pi.
\]
Substituting (\ref{DM_parameter}) gives 
the corresponding twistor line parametrised by $[\pi]\in\PP^1$ 
\be
\label{sl3curves}
p^{\alpha}=[(z\cdot \pi)x^A-\pi^A, z\cdot \pi], \quad l_\alpha=[\pi_A, -\pi\cdot x],
\ee
where the spinor indices are raised and lowered with $\epsilon^{AB}$ and its inverse, and  $z\cdot x\equiv z_Ax^A$. 

We shall now derive the expression for the conformal structure. According
to the Nonlinear Graviton prescription of Penrose \cite{penrose} a vector
$V\in \Gamma(T_mM)$ is null if the corresponding section
of the normal bundle $N(L_m)={\mathcal O}(1)\oplus{\mathcal O}(1)$ has a 
single zero.  To compute the normal bundle, let $([l(\pi, P, L)], 
[p(\pi, P, L)])$
be the twistor line corresponding to a point $m=(P, L)$ in $M$. The neighbouring line is $([l+\delta l], [p+\delta p])$, where
\[
\delta l=\delta P\wedge \pi, \quad
\delta p= (\delta L\cdot \pi)P+(L\cdot\pi) \delta P-\delta (L\cdot P)\pi.
\]
The lines  $(l+\delta l, p+\delta p)$  and $(l, p)$ intersect at one point
which correspond to some particular value of $\pi$. Therefore
\[
l+\delta l\sim l, \quad\mbox{so}\quad \pi\sim\delta P=[\delta x^1, \delta x^2, 0].
\]
The other condition is $p+\delta p\sim p$ which holds iff 
\[
0=p\wedge \delta p=(L\cdot \pi)^2P\wedge \delta P-(L\cdot P)
(\delta L\cdot \pi)\pi\wedge P-(L\cdot \pi)\delta (L\cdot P)P\wedge \pi-
(L\cdot P)(L\cdot \pi) \pi\wedge \delta P.
\]
Substituting $\pi\sim\delta P$, we find that all terms on the RHS are proportional to $P\wedge \delta P=[0, 0, x\cdot dx]$.
Therefore
\[
(L\cdot \pi)^2-(L\cdot \pi)\delta(L\cdot P)+(L\cdot P)(\delta L\cdot \pi)=0,
\]
together with $L\cdot P=1$.  This gives the conformal structure 
(\ref{dm_metric}).
%%%%%%%%%%%%%%%%%%%%%%%%%%%%%%%%%%%%%%%%%%%%%%%%%%%%%%%%%%%%%%
\subsection{Hyper--Hermitian structure}
A pseudo--hyper--complex structure on a four manifold $M$ is a triple of endomorphisms
$I, S, T$ of $TM$ which satisfy
\[
I^2=-Id, \quad S^2=T^2=Id, \quad IST=Id,
\]
and such that $aI+bS+cT$ is an integrable complex structure for any point 
on the hyperboloid $a^2-b^2-c^2=1$.
A neutral signature metric $g$ on a pseudo--hyper--complex four--manifold is pseudo-hyper-Hermitian
if
\[
g(V, V)=g(IV, IV)=-g(SV, SV)=-g(TV, TV)
\]
for any vector field $V$ on $M$. There is a unique conformal structure compatible
with each pseudo--hyper--complex structure. With a natural choice of orientation
which makes the fundamental two--forms of $I, S, T$  self--dual, 
this conformal structure is anti--self--dual.
\begin{prop}
\label{propHH}
The Einstein metric (\ref{special_ein}) on $SL(3)/GL(2)$ is pseudo--hyper--Hermitian.
\end{prop}
\noindent
{\bf Proof.}
First a note about conventions.
Our model metric is ASD, so it is the primed Weyl spinor which vanishes. This makes all indices on both the EW space and the projective surface primed. Lets therefore swap the role of primed and unprimed indices. The null frame for the 4-metric is
\[
e^{0'A}=dx^A, \quad e^{1'A}=dz^A+z^A(z\cdot dx), \quad\mbox{so that}\quad
g=\epsilon_{A'B'}\epsilon_{AB}e^{A'A}e^{B'B}.
\]
Thus the forms $\Sigma=dx^0\wedge dx^1$ and $\Omega=dz_A\wedge dx^A$ are ASD. The basis of SD two forms is spanned by
\[
dx\wedge dq+ q^2 dx\wedge dy,\quad
dx\wedge dp-dy\wedge dq+2 pq dx\wedge dy,
\quad
-dy\wedge dp+ p^2 dx\wedge dy
\]
or, in a more compact notation, by
$\Sigma^{AB}=dx^{(A}\wedge dz^{B)}+z^Az^B \Sigma$.
We can verify that
\be
\label{lie_form}
d\Sigma^{AB}+2{{\mathcal{A}}}\wedge\Sigma^{AB}=0,
\ee
where ${{\mathcal A}}= z_Adx^A$ is such that $d{{\mathcal{A}}}= \Omega$. 
The condition (\ref{lie_form}) is neccesary and sufficient for (para) hyper--Hermiticity
\cite{boyer, D99}.
Thus the ASD
Maxwell fields arising from the para--K\"ahler structure on $M$, and the para hyper--Hermitian structure coincide.
To this end note that the twistor distribution form $(M,  g)$ is
\be
\label{tdistribution}
L_{0'}=\pi\cdot\frac{\p}{\p x}+(z\cdot\pi) z\cdot\frac{\p}{\p z}, \quad
L_{1'}=\pi\cdot\frac{\p}{\p z}.
\ee
It is Frobenius integrable, as $[L_{0'}, L_{1'}]=-(\pi\cdot z)L_{1'}$. It also does not contain the vertical $\p/\partial \pi$ terms which again confirms the hyper--Hermiticity
of $(M, g)$ (see Lemma 2 in \cite{D99} and Theorem  7.1 in \cite{Cal2}). The SD part of the
spin connection is given in terms of ${\mathcal A}$ as
$\Gamma_{A'ABC}=-2{\mathcal A}_{A'(B}\epsilon_{C)A}$.
\koniec
In the next section we shall show how to encode
${\mathcal A}$ in the twisted--photon Ward bundle over the twistor space
of $(M, g)$.
\subsection{The twisted photon}
The twistor space $F_{12}$ described in \S\ref{twist_SSS} is the projectivised tangent bundle ${T(\PP^2)}^*$ of the minitwistor space of the flat projective structure: a point in $F_{12}$ consists of $l\in\PP^2$, and a direction through $l$. Thus the twistor space of $M$ is the correspondence space
of $\PP^2$ and ${\PP^2}^*$. 
There are many open sets needed to cover
$\PP(T\PP^2)$, but it is sufficient to consider two:
$U$, where $(l_1, \neq 0, p^2\neq 0)$, and $(l_2/l_1, l_3/l_1, p^3/p^2)$ are coordinates, and $\widetilde{U}$ where
$(l_1\neq 0, p^3\neq 0)$, and  $(l_2/l_1, l_3/l_1, p^2/p^3)$
are coordinates. Now consider the total
space of $T\PP^2$ (or perhaps it is $T\PP^2$ tensored
with some power of the canonical bundle to make it trivial on twistor
lines), and restrict it to the intersection of (pre--images in
$T\PP^2$
of) $U$ and $\widetilde{U}$. The coordinates on $T\PP^2$ in these
region are $(l_2/l_1, l_3/l_1, p^2/p^1, p^3/p^1)$, and the fiber
coordinates $\tau$ over $U$ and $\tilde{\tau}$ over 
$\widetilde{U}$ are related by\footnote{Here we are following Ward \cite{wardtf},
and thinking of a $\C^*$ bundle.}
\[
\tilde \tau=\exp(F)\tau, \quad\mbox{where}\quad  
F=\ln{(p_2/p_3)}.
\]
Now we follow the procedure of \cite{wardtf}: restrict $F$ to a twistor line,
and split it.
The holomorphic splitting is $F=H-\widetilde{H}$, where
$H=\ln{(p_2)}$ is holomorphic in the pre--image of $U$ in the correspondence space, and 
$\widetilde{H}=\ln{(p_3)}$ is holomorphic in the pre--image of
$\widetilde{U}$. Note that $F$ is a twistor  function, but 
$H, \widetilde{H}$ are not. Therefore
$L_{A'}F=0$, where the twistor distribution $L_{A'}$
is given by (\ref{tdistribution}). This, 
together with the Liouville theorem
implies
that
\[
L_{A'}H=L_{A'}\widetilde{H}=\pi^A\mathcal{A}_{A'A}
\]
for some one--form $\mathcal{A}$ on $M$, 
as the LHS is holomorphic on  $\CP^1$ and homogeneous of degree
one. To construct this one--form recall the parametrisation
of twistor curves (\ref{sl3curves}). This gives
\[
H=\ln{(z\cdot\pi)}, \quad\widetilde{H}=\ln{((z\cdot\pi)x^1-\pi^1)}
\]
and
\[
L_{1'}(H)=L_{1'}(\widetilde{H})=0, \quad
L_{0'}(H)=L_{0'}(\widetilde{H})=\pi\cdot z.
\]
Therefore ${\mathcal A}_{1'A}=0, {\mathcal A}_{0'A}=z_A$
which gives ${\mathcal A}=z_Adx^A$, and $d{\mathcal A}$
is indeed the ASD para--K\"ahler structure.
%which I THINK is 
%\be
%\label{ptp}
%p=[p^1, p^2, p^3], \quad Z=-\frac{p^3}{p^2}=\frac{z^0+X %z^1}{-1+(z^0+X z^1)x^1}.
%\ee
\subsection{Factoring $SL(3)/GL(2)$ to Einstein-Weyl}
\label{neat2}
If a metric with ASD Weyl tensor has more than one conformal symmetry, then distinct Einstein--Weyl structures are obtained on the space of orbits of conformal Killing vectors which are not conjugate with respect to an isometry \cite{PT}. We can thus classify the Einstein--Weyl structures obtainable from $SL(3)/GL(2)$ by first classifying its symmetries up to conjugation.
\begin{prop}
The non--trivial Einstein--Weyl structures obtainable from $SL(3)/GL(2)$ by the Jones--Tod correspondence consist of a two--parameter family, and two additional cases which do not belong to this family.
\end{prop}
\noindent
{\bf Proof. }Since we have an isomorphism between the Lie algebra of projective vector fields on $(N,[\nabla])$ and the Lie algebra of Killing vectors on $(M,g)$, the problem of classifying the symmetries of $M=SL(3)/GL(2)$ is reduced to a classification of the infinitesimal  projective symmetries of $\RP^2$, i.e. the near--identity elements of $SL(3)$, up to conjugation. Non--singular complex matrices are determined up to similarity by their Jordan normal form (JNF). While real matrices do not have such a canonical form, all of the information they contain is determined (up to similarity) by the JNF that they would have if they were considered as complex matrices. Thus we can still discuss the JNF of a real matrix, even if it cannot always be obtained from the real matrix by a real similarity transformation. The possible non--trivial Jordan normal forms of matrices in $SL(3)$ are shown below.
\[
\begin{pmatrix}\lambda & 0 & 0\\
0 & \mu & 0\\
0 & 0 & 1/\lambda\mu
\end{pmatrix}
\quad
\begin{pmatrix}\lambda & 0 & 0\\
0 & \lambda & 0\\
0 & 0 & 1/\lambda^2
\end{pmatrix}
\quad
\begin{pmatrix}\lambda & 1 & 0\\
0 & \lambda & 0\\
0 & 0 & 1/\lambda^2
\end{pmatrix}
\quad
\begin{pmatrix}1 & 1 & 0\\
0 & 1 & 0\\
0 & 0 & 1
\end{pmatrix}
\quad
\begin{pmatrix}1 & 1 & 0\\
0 & 1 & 1\\
0 & 0 & 1
\end{pmatrix}
\]

It is possible that two matrices in $SL(3)$ with the same JNF may be related by a complex similarity transformation, and thus not conjugate in $SL(3)$. However, if the JNF is a real matrix, then the required similarity transformation just consists of the eigenvectors and generalised eigenvectors of the matrix, which must also be real since they are defined by real linear simultaneous equations. This means we only have to worry about matrices with complex eigenvalues, and since these occur in complex conjugate pairs, they will only be a problem when we have three distinct eigenvalues.

In this case, we can always make a real similarity transformation such that the matrix is block diagonal, with the real eigenvalue in the bottom right. Then we have limited choice from the $2\times 2$ matrix in the top left. Let us parametrise such a $2\times 2$ matrix by $a,\,b,\,c,\,d\in\mathbb{R}$ as follows:
\[
\quad
\begin{pmatrix}1+a\epsilon & b\epsilon \\
c\epsilon & 1+d\epsilon 
\end{pmatrix}.
\]
This has characteristic polynomial
\[
\chi(\lambda)=\lambda^2-(2+\epsilon(a+d))\lambda+1+(a+d)\epsilon+(ad-bc)\epsilon^2.
\]
Evidently the important degrees of freedom are $a+d$ and $ad-bc$, so we can use these to encode every near--identity element of the class with three distinct eigenvalues. The bottom--right entry will be determined by our choice of $a+d$ and $ad-bc$.

Taking a projective vector field on $\mathbb{RP}^2$, we can find the corresponding Killing vector on $SL(3)/GL(2)$ using (\ref{symlift}), and factor to Einstein--Weyl using (\ref{EWgen}). We find by explicit calculation that vector fields arising from the second and fourth JNFs above give trivial Einstein--Weyl structures, so restricting to the non--trivial cases we have a two--parameter family of Einstein--Weyl structures coming from the first class, and two additional Einstein--Weyl structures coming from the third and fifth, as claimed.
\begin{flushright}
$\square$
\par\end{flushright}
\subsection{Mini--twistor correspondence}
\label{mini_twistor}
Below we investigate a one--parameter subfamily of the two--parameter family. We use the holomorphic vector field on 
the twistor space
$F_{12}$  (see \S\ref{twist_SSS})
corresponding to the chosen symmetry, and reconstruct the conformal structure $[h]$ on $N$ using minitwistor curves 
(in the sense of \cite{hitchin})
on the space of orbits. Take $a\in \mathbb{R}$ and
\be
\label{modelK}
K=P^1\frac{\partial}{\partial P^1} - L_1\frac{\partial}{\partial L_1}+aP^2\frac{\partial}{\partial P^2} - aL_2\frac{\partial}{\partial L_2},
\ee
%corresponding to the matrix
%\[
%M=
%\begin{pmatrix}1+\epsilon & 0 & 0\\
%0 & 1+a\epsilon & 0\\
%0 & 0 & 1
%\end{pmatrix}.
%\]
%Note that we have chosen inhomogeneous coordinates, so $M$ need not be in $SL(3)$.
In order to preserve the relations
\[
p\cdot L=0,\quad P\cdot l=0,\quad p\cdot l=0, 
\]
the corresponding holomorphic action on $(p,l)$ must be $p\mapsto Mp$, $l\mapsto M^{-1}l$, thus the holomorphic vector field $\mathcal{K}$ on $F_{12}$ is
\[
\mathcal{K}=p^1\frac{\partial}{\partial p^1} - l_1\frac{\partial}{\partial l_1}+ap^2\frac{\partial}{\partial p^2} - al_2\frac{\partial}{\partial l_2}.
\]

In order to factor $F_{12}$ by this vector field, we must find invariant minitwistor coordinates $(Q,R)$. In addition to satisfying $\mathcal{K}(Q)=\mathcal{K}(R)=0$, they must be homogeneous of degree zero in $(P,L)$. We choose
\[
Q=\frac{p^1l_1}{p^2l_2},\quad R=\frac{(l_1)^a}{l_2(l_3)^{a-1}}.
\]
Substituting in our parametrisation (\ref{sl3curves}) and using the freedom to perform a Mobius transformation on $\pi$, we obtain
\begin{align}
\label{QR}
Q &= \frac{(\lambda t-u-1)\lambda}{v\lambda+\lambda-\frac{uv}{t}}\\
R &= \lambda^a\Big(-\lambda-\frac{v}{t}\Big)^{1-a},\nonumber
\end{align}
where we have defined $\lambda=\pi_0/\pi_1$, and the Einstein-Weyl coordinates
\[
u=xp,\quad v=yq, \quad t=x^aq.
\]
Note these are invariants of the Killing vector (\ref{modelK}).

Next we wish to use these minitwistor curves to reconstruct the conformal structure of the Einstein-Weyl space. In doing so we follow \cite{PT}. The tangent vector field to a fixed curve is given by
\[
T=\frac{\p Q}{\p \lambda} \frac{\p}{\p Q} + \frac{\p R}{\p \lambda} \frac{\p}{\p R},
\]
Hence we can write the normal vector field as
\begin{align*}
N &=dQ\frac{\p}{\p Q} + dR\frac{\p}{\p R} \enskip \mathrm{mod}\, T\\
&= \bigg(\frac{\p R}{\p \lambda}\bigg)^{-1}\bigg(dQ\frac{\p R}{\p \lambda}-dR\frac{\p Q}{\p \lambda}\bigg)\frac{\p}{\p Q},
\end{align*}
where
\[
dQ=\frac{\p Q}{\p u}du + \frac{\p Q}{\p v}dv + \frac{\p Q}{\p t}dt
\]
and similarly for $dR$. Calculating $N$ using (\ref{QR}), we find
\[
N\propto(A\lambda^2+B\lambda+C)\frac{\p}{\p Q},
\]
where
\begin{align*}
A &= t^2(v+1)dt-t^3dv, \\
B &= -2tuvdt + t^2(a+2u)dv -t^2du, \\
C &= uv(1+u)dt - tu(1+u)dv - atvdu.
\end{align*}
The discriminant of this quadratic in $\lambda$ then gives a representative $h\in[h]$ of our conformal structure:
\be
\begin{split}
h=4(u^2v+uv^2+uv)dt^2 - 4tv(a(v+1)+u)dtdu + 4tu(u-av+2v+1)dtdv \\
-t^2du^2 + 2t^2(2av+a+2u)dvdu - t^2(a^2+4u(a-1))dv^2.
\end{split}
\ee
This is the same conformal structure that we obtain by Jones-Tod factorisation of $SL(3)/GL(2)$ by (\ref{modelK}) using the formula (\ref{EWgen}).
\section{Einstein--Weyl from a third order ODE}
\label{ODEs}
Let 
\be
\label{ode1}
Y''=F(X, Y, Y')
\ee
be the second order ODE dual to the ODE (\ref{odealice}) defining the projective
structure on $N$. The integral curves of (\ref{ode1}) are twistor curves corresponding with normal bundle ${\mathcal O}(1)$ corresponding to points in $N$. The construction
of \cite{DM} suggests an implict way of extending the ODE (\ref{ode1}) to a system
\be
\label{ode2}
Y''=F(X, Y, Y'), \quad Z''=G(X, Y, Z, Y', Z')
\ee
with vanishing Wilczynski invariants \cite{Grossman, Doubrov, CDT13} ,  such that the solution space
to (\ref{ode2}) carries an ASD structure. Assume that
this ASD Einstein structure admits a Killing vector $K$. This corresponds to
a Lie point symmetry of the system (\ref{ode2}) \cite{CDT13}, and 
(by Theorem 3.10
of \cite{DM}) to a Lie point symmetry of (\ref{ode1}). Without loss of generality
we can now choose the coordinates $(X, Y)$ such that $F=F(Y, Y')$ in 
(\ref{ode1}). We shall now construct a third order ODE which, by a combination of the Cartan
correspondence \cite{Cartan} and the Jones--Tod construction \cite{JT} 
gives rise to an Einstein--Weyl structure on the space of orbits of $K$.
The third order ODE is obtained by differentiating the second equation
in (\ref{ode2}) with respect to $X$, eliminating $Y''$ using the first equation,  and
finally eliminating $Y'$ using the second equation.  
This yields\footnote{We are taking the point symmetry of (\ref{ode1}) to be $\p/\p X$. 
which lifts to a point symmetry of (\ref{ode2}) of the form
$\p/\p X+J(X, Y, Z)\p/\p Z$, so that there can be an explicit $X$--dependence
in $G$ and therefore in $H$.}
\be
\label{3rdorder}
Z'''=H(X, Z, Z', Z'')\quad\mbox{where}\quad H=\Big(\frac{\p G}{\p X}
+\frac{\p G}{\p Y}Y'+\frac{\p G}{\p Y'}F+\frac{\p G}{\p Z}Z'+
\frac{\p G}{\p Z'}G\Big)|_{Y'=Y'(X, Z, Z', Z'')}.
\ee
\subsection{Example 1}
Consider the system of ODEs \cite{CDT13}
\be
\label{sl3ode}
Y''=0, \quad Z''=\frac{2(Z')^2Y'}{ZY'-1}
\ee
It  admits $SL(3)$ as its group of point symmetries, and its solution
space carries the the conformal structure (\ref{dm1}). The procedure
described above leads to the Schwartzian third order  ODE
\[
Z'''=\frac{3 (Z'')^2}{2Z'}.
\]
The integral curves are of the form
\[
Z=\frac{aX+b}{cX+d}, \quad{where}\quad ad-bc=1
\]
and the corresponding Einstein--Weyl structure is flat.

To this end we show how to derive the system (\ref{sl3ode}) from the twistor curves
(\ref{sl3curves}). Let $X=\pi_1/\pi_0$ be the affine coordinate
on the $\PP^1$ parametrising (\ref{sl3curves}). Then
\[
l=[1,X, Y]=[1, X, -x^0-Xx^1], \quad\mbox{so that}\quad Y''=0.
\]
To recover the second ODE in (\ref{sl3ode}) we need an explicit identification
of $F_{12}$ with $\PP(T\PP^2)$, which is 
\[
p=[p^1, p^2, p^3], \quad Z=-\frac{p^3}{p^2}=\frac{z^0+X z^1}{-1+(z^0+X z^1)x^1}.
\]
The pair $(Y(X), Z(X))$ is the general solution to (\ref{sl3ode}) so this ODE is
indeed the {\it right one} for the twistor curves in $F_{12}$.
\subsection{Example 2.}
Another example is  given by the system \cite{CDT13}
\[
Y''=0, \quad Z''=-(Z'+\sqrt{(Y')^2-1})^2.
\]
The resulting ASD is not Einstein,  and so does not 
fit into the scheme of \cite{DM} (it does however 
admit an ASD $\beta$--distribution coming from a null 
Killing vector).  The corresponding third order ODE (\ref{3rdorder})
is
\[
Z'''=2(Z'')^{3/2},
\]
and the resulting Einstein--Weyl structure is
the Lorentzian Nil structure on the Heisenberg group \cite{T}.

%{\color{blue} There are two major  missing points in all this
%\begin{enumerate}
%\item Given $F$ in (\ref{ode1}) find $G$ in (\ref{ode2}) such that
%the resulting ASD structure is that of \cite{DM}.
%\item Show that the Wilczynski condtions for the system (\ref{ode2}) imply
%the Wunschman and the Cartan condition for the 3rd order ODE (\ref{3rdorder}).
%\end{enumerate}}

\section{Kaluza-Klein lift to an $S^{1}$-bundle over $M$}
\label{KKsection}

In this section we consider the generalisation of the four-manifold
$M$ to any even dimension $2n$, starting with a projective structure
of dimension $n$. This is discussed in the appendix of \cite{DM}. The construction
of \cite{DM} offers two different perspectives on $(M, g)$: one in terms of a gauge theory on $(N,[\nabla])$ and the other as a quotient of the Cartan bundle of $(N, [\nabla])$. It is the latter approach that we will take here. The aim is to construct a principal circle bundle over the $2n$--dimensional $M$ with a metric on the $2n+1$--dimensional total space which we will prove to be Einstein. This metric is interesting in the context of Kaluza-Klein theory.

\subsection{The Cartan geometry of a projective structure}
Our description follows \cite{Bryant} and \cite{DM}.
A projective structure in dimension $n$ admits a description as a Cartan geometry by virtue of being modelled on $\RP^n$, since this can be viewed as a homogeneous space. Let $H$ be the stabiliser subgroup of the point $[1,0,\dots,0]$
under the action of $SL(n+1,\mathbb{R})$ on homogeneous coordinates
$[X^{0},\dots,X^{n}]$ in $\mathbb{RP}^{n}$ via its fundamental representation. Then by the orbit stabiliser theorem $\RP^n=SL(n+1,\R)/H$. The elements of $H$
take the general form
\[
\begin{pmatrix}\mathrm{det}(a)^{-1} & b\\
0 & a
\end{pmatrix}
\]
for any $a\in GL(n,\mathbb{R})$ and $b\in\mathbb{R}_{n}$.

\begin{defi}
The Cartan geometry $(\pi:P\mapsto N,\theta)$ of a projective structure $(N,[\nabla])$ in dimension $n$ consists of a principal right $H$--bundle over $N$ carrying an $\mathfrak{sl}(n+1,\mathbb{R})$--valued one--form $\theta$ with the following properties:
\begin{enumerate}
\item{$\theta(X_v)=v$ for fundamental vector fields $X_v$ on $P$ corresponding to $v\in\mathfrak{h}$;}
\item{$\theta_u:T_uP\rightarrow\mathfrak{sl}(n+1,\R)$ is an isomorphism for all $u\in P$;}
\item{$R^{*}_h\theta=\mathrm{Ad}(h^{-1})\theta=h^{-1}\theta h$ for all $h\in H$.}
\end{enumerate}
The form $\theta$ is called the Cartan connection.
\end{defi}

\noindent We write the Cartan connection as
\[
\theta=\begin{pmatrix}-\mathrm{tr}(\phi) & \eta\\
\omega & \phi
\end{pmatrix}
\]
with the one-forms $\phi$, $\eta$ and $\omega$ taking values
in $\mathfrak{gl}(n,\mathbb{R})$, $\mathbb{R}^{n}$ and $\mathbb{R}_{n}$
respectively. The curvature two-form takes the general form
\begin{equation}
\Theta={d}\theta+\theta\wedge\theta=\begin{pmatrix}0 & L(\omega\wedge\omega)\\
0 & W(\omega\wedge\omega)
\end{pmatrix},\label{eq:curvature_2-form}
\end{equation}
where $L$ and $W$ are curvature functions valued in $\mathrm{Hom}(\R^n\wedge\R^n,\R_n)$ and $\mathrm{Hom}(\mathbb{R}^{n}\wedge\mathbb{R}^{n},\mathbb{R}_{n}\otimes\mathbb{R}^{n})$ respectively.

\subsection{Construction of $(M,g)$ as quotient of the Cartan bundle}
We can embed $GL(n,\mathbb{R})$ as a subgroup of $H$ via the map
\begin{equation}
GL(n,\mathbb{R})\ni a\longmapsto\begin{pmatrix}\mathrm{det}(a)^{-1} & 0\\
0 & a
\end{pmatrix}\in H,\label{eq:gl(n)_embed_H}
\end{equation}
finding that its adjoint action on $\theta$ preserves the natural
contraction $\eta\omega:=\sum_{i}\eta_{(i)}\otimes\omega^{(i)}$,
where $\eta_{(i)}$ and $\omega^{(i)}$ are the one-forms denoting
the components of $\eta$ and $\omega$ in $\mathfrak{sl}(n+1,\mathbb{R})$. We deduce that the quotient $q:P\mapsto P/GL(n,\mathbb{R})$
carries a natural metric $g$ of signature $(n,n)$ and a two-form
$\Omega$ which are such that
\begin{eqnarray*}
q^{*}g & = & \mathrm{Sym}(\eta\omega)\\
q^{*}\Omega & = & \mathrm{Ant}(\eta\omega),
\end{eqnarray*}
where $\mathrm{Sym}$ and $\mathrm{Ant}$ denote the symmetric and anti--symmetric parts of the $(0,2)$ tensor $\eta\omega$.
It is shown in \cite{DM} that we can choose local coordinates such that they take the forms
(\ref{g}) and (\ref{Omega}) respectively,
where the indices $A,\,B$ now run from $0$ to $n-1$. The pair $(g,\Omega)$ forms a so--called bi--Lagrangian structure on $M$. The fact that $\Omega$ is symplectic follows from the Bianchi identity.

In fact, the metric and symplectic form turn out to belong to a one--parameter family, which can be written in local coordinates as
\begin{eqnarray}
\label{g_Lambda}
g_\Lambda &=& dz_A\odot dx^A-(\Gamma^C_{AB}z_C-\Lambda z_Az_B-\Lambda^{-1}P_{AB})dx^A\odot dx^B \\
\Omega_\Lambda &=& dz_A\odot dx^A + \frac{1}{\Lambda}P_{AB}dx^A\wedge dx^B,\quad A,B=0,\dots,n-1 \nonumber
\end{eqnarray}
for $\Lambda\neq 0$. Metrics of the form (\ref{g_Lambda}) are a subclass of so--called Osserman metrics. More details can be found in \cite{osserman}. They are all Einstein with non--zero scalar curvature (calculated in  Appendix A), but for $\Lambda\neq 1$ the relation to projective geometry is lost. For the remainder of this section we will consider the full family $\{g_\Lambda\}$, although projective invariance will only hold in the case $\Lambda=1$.

\subsection{Flat case}

We first note that taking the flat projective structure on $\mathbb{RP}^{n}$
results in a one--parameter family of $2n$-dimensional Einstein spaces $M$ which are Kaluza-Klein
reductions of quadrics in $\mathbb{R}^{n+1,n+1}$. The metric and symplectic form on $M$ reduce to 
\begin{eqnarray*}
g_{\Lambda} & = & {d}x^{A}\odot {d}z_{A}+\Lambda z_{A}z_B{d}x^{A}\odot dx^B\\
\Omega_{\Lambda} & = & {d}z_{A}\wedge{d}x^{A},\qquad A,B=0,\dots,n-1.
\end{eqnarray*}

\begin{prop}The Einstein spaces $M$ corresponding to $\mathbb{RP}^{n}$
are projections from the $2n+1$-dimensional quadrics $Q\subset\mathbb{R}^{n+1,n+1}$
given by $X^{\alpha}Y_{\alpha}=\frac{1}{\Lambda}$, where $X,Y\in\mathbb{R}^{n+1}$
are coordinates on $\mathbb{R}^{n+1,n+1}$ such that the metric is
given by 
\[
G={d}X^{\alpha}{d}Y_{\alpha},
\]
 under the embedding
\begin{eqnarray}
X^{\alpha} & = & \begin{cases}
x^{A}\mathrm{e}^{\tau}, & \alpha=A=0,\dots,n-1\\
\mathrm{e}^{\tau}, & \alpha=n
\end{cases}\nonumber \\
Y_{\alpha} & = & \begin{cases}
z_{A}\mathrm{e}^{-\tau}, & \alpha=A=0,\dots,n-1\\
\mathrm{e}^{-\tau}\Bigl(\frac{1}{\Lambda}-x^{C}z_{C}\Bigr), & \alpha=n
\end{cases}\label{eq:embedding-1}
\end{eqnarray}
following Kaluza-Klein reduction by the vector $\frac{\partial}{\partial\tau}.$
\end{prop}
\noindent
\textbf{Proof. }We find the basis of coordinate 1-forms $\{{d}X^{\alpha},{d}Y_{\alpha}\}$
to be
\begin{eqnarray*}
{d}X^{\alpha} & = & \begin{cases}
\mathrm{e}^{\tau}({d}x^{A}+x^{A}{d}\tau), & \alpha=A=0,\dots,n-1\\
\mathrm{e}^{\tau}{d}\tau, & \alpha=n
\end{cases}\\
{d}Y_{\alpha} & = & \begin{cases}
\mathrm{e}^{-\tau}({d}z_{A}-z_{A}{d}\tau), & \alpha=A=0,\dots,n-1\\
-\mathrm{e}^{-\tau}\biggl[\Bigl(\frac{1}{\Lambda}-x^{C}z_{C}\Bigr){d}\tau+x^{C}{d}z_{C}+z_{C}{d}x^{C}\biggr], & \alpha=n.
\end{cases}
\end{eqnarray*}
The metric is then given by
\begin{eqnarray*}
G & = & \mathrm{e}^{\tau}({d}x^{A}+x^{A}{d}\tau)\mathrm{e}^{-\tau}({d}z_{A}-z_{A}{d}\tau)\ -\ \mathrm{e}^{\tau}{d}\tau\mathrm{e}^{-\tau}\biggl[\Bigl(\frac{1}{\Lambda}-x^{C}z_{C}\Bigr){d}\tau+x^{C}{d}z_{C}+z_{C}{d}x^{C}\biggr]\\
 & = & {d}x^{A}{d}z_{A}\ +\ (x^{A}{d}z_{A}-z_{A}{d}x^{A}){d}\tau\ -\ (x^{A}z_{A}){d}\tau^{2}\ -\ {d}\tau\biggl[\Bigl(\frac{1}{\Lambda}-x^{C}z_{C}\Bigr){d}\tau+x^{C}{d}z_{C}+z_{C}{d}x^{C}\biggr]\\
 & = & {d}x^{A}{d}z_{A}\ -\ \frac{1}{\Lambda}{d}\tau^{2}\ -\ 2z_{A}{d}x^{A}{d}\tau\\
 & = & {d}x^{A}{d}z_{A}\ +\ \Lambda(z_{A}{d}x^{A})^{2}\ -\ \Lambda\Bigl(\frac{{d}\tau}{\Lambda}+z_{A}{d}x^{A}\Bigr)^{2},
\end{eqnarray*}
which is clearly going to give $g_{\Lambda}$ under Kaluza-Klein reduction
by $\frac{\partial}{\partial\tau}$.
\begin{flushright}
$\square$
\par\end{flushright}

Note that the symplectic form $\Omega$ is the exterior derivative
of the potential term $z_{A}{d}x^{A}$, implying a possible
generalisation to the curved case.

\subsection{Curved case}\label{kk_curved_case}

We now return to a general projective structure $(N,[\nabla])$. Since the
symplectic form picks out the antisymmetric part of the Schouten tensor,
it has the fairly simple form
\[
\Omega_{\Lambda}={d}z_{A}\wedge{d}x^{A}-\frac{\partial_{[A}\Gamma_{B]C}^{C}}{\Lambda(n+1)}{d}x^{A}\wedge{d}x^{B}.
\]
By inspection, this can be written $\Omega_{\Lambda}={d}\mathcal{A}$,
where
\[
\mathcal{A}=z_{A}{d}x^{A}-\frac{\Gamma_{AC}^{C}}{\Lambda(n+1)}{d}x^{A}.
\]
This is a trivialisation of the Kaluza-Klein bundle which we are about
to construct. Note that for $\Lambda=1$, $\Omega$ and
$\mathcal{A}$ remain unchanged under a change of projective connection (\ref{proj_change}).

Motivated by the Kaluza-Klein reduction in the flat case, we consider
the following metric.
\begin{theo}
\label{theokk}
Let $g_{\Lambda}$ be the Einstein metric
(\ref{g_Lambda}) corresponding to the projective structure $(N, [\nabla])$.
The metric
\begin{equation}
\mathcal{G}_{\Lambda}=g_{\Lambda}-\Lambda\Bigl(\frac{{d}t}{\Lambda}+\mathcal{A}\Bigr)^{2}\label{eq:G}
\end{equation}
on a principal circle bundle $\sigma:\mathcal{Q}\mapsto M$ is Einstein,
with Ricci scalar $2n(2n+1)\Lambda$.
\end{theo}
\noindent
\textbf{Proof.} We prove this using the Cartan formalism. Our treatment
parallels a calculation by Kobayashi \cite{Kob}, who considered
principal circle bundles over K\"ahler manifolds in order to study the
topology of the base. Note that we temporarily surpress the constant
$\Lambda$, writing $\mathcal{G}\equiv\mathcal{G}_{\Lambda}$ and
$g\equiv g_{\Lambda}$, since the proof applies to any choice $\Lambda\neq0$
within this family. Consider a frame
\begin{equation}
e^{a}=\begin{cases}
{d}x^{A}, & a=A=0,\dots,n-1\\
{d}z_{A}-(\Gamma_{AB}^{C}z_{C}-\Lambda z_{A}z_{B}-\Lambda^{-1}P_{AB}){d}x^{B}, & a=A+n=n,\dots,2n-1,
\end{cases}\label{eq:basis}
\end{equation}
and in this basis the metric takes the form
\begin{equation}
g=e^{0}\odot e^{n}+\dots+e^{n-1}\odot e^{2n-1}.\label{eq:g_cov_const}
\end{equation}
We are interested in the metric
\[
\mathcal{G}=g-e^{t}\odot e^{t},
\]
where
\[
e^{t}=\sqrt{\Lambda}\biggl(\frac{{d}t}{\Lambda}+{\mathcal A}\biggr).
\]
 We reserve Roman indices $a,b,\dots$ for the $2n$-metric components
$0,\dots,2n-1$ and allow greek indices $\mu,\nu,\dots$ to take values
$0,1,\dots,2n$. The dual basis to $\{e^{\mu}\}$ will be denoted
$\{E_{\mu}\}$ and will act on functions as vector fields in the usual
way. We wish to find the new connection 1-forms $\hat{\psi}_{\ \nu}^{\mu}$
(defined by ${d}e^{\mu}=-\hat{\psi}_{\ \nu}^{\mu}\wedge e^{\nu}$)
in terms of the old ones $\psi_{\ b}^{a}$ (defined by ${d}e^{a}=-\psi_{\ b}^{a}\wedge e^{b}$).
Hence we examine%
\footnote{Note that our conventions are $({d}\omega)_{ab\dots c}=\partial_{[a}\omega_{b\dots c]},$
$(\eta\wedge\omega)_{a\dots d}=\eta_{[a\dots b}\omega_{c\dots d]},$
$\omega=\omega_{a\dots b}{d}x^{a}\wedge\dots\wedge{d}x^{b},$
and $F_{ab}{d}x^{a}\wedge{d}x^{b}=F_{[ab]}{d}x^{a}\otimes{d}x^{b}$
implying ${d}x^{a}\wedge{d}x^{b}=\frac{1}{2}({d}x^{a}\otimes{d}x^{b}-{d}x^{b}\otimes{d}x^{a})$.%
} ${d}e^{t}$ to find $\hat{\psi}_{\ a}^{t}.$
\[
{d}e^{t}=\sqrt{\Lambda}{d}A=\sqrt{\Lambda}\Omega_{ab}e^{a}\wedge e^{b}=-\hat{\psi}_{\ a}^{t}\wedge e^{a}\quad\implies\quad\hat{\psi}_{\ a}^{t}=\sqrt{\Lambda}\Omega_{[ab]}e^{b}=\sqrt{\Lambda}\Omega_{ab}e^{b},\qquad\hat{\psi}_{\ t}^{a}=\sqrt{\Lambda}\Omega_{\ b}^{a}e^{b}.
\]
Since ${d}e^{a}$ is unchanged, we have that
\[
\hat{\psi}_{\ t}^{a}\wedge e^{t}+\hat{\psi}_{\ b}^{a}\wedge e^{b}=\psi_{\ b}^{a}\wedge e^{b},
\]
thus
\[
\hat{\psi}_{\ b}^{a}\wedge e^{b}=\psi_{\ b}^{a}\wedge e^{b}-\sqrt{\Lambda}\Omega_{\ b}^{a}e^{b}\wedge e^{t}\qquad\implies\qquad\hat{\psi}_{\ b}^{a}=\psi_{\ b}^{a}+\sqrt{\Lambda}\Omega_{\ b}^{a}e^{t}.
\]

We now calculate the curvature 2-forms $\hat{\Psi}_{\ \nu}^{\mu}={d}\hat{\psi}_{\ \nu}^{\mu}+\hat{\psi}_{\ \rho}^{\mu}\wedge\hat{\psi}_{\ \nu}^{\rho}=\frac{1}{2}\mathcal{R}_{\rho\sigma\nu}^{\ \ \ \ \mu}e^{\rho}\wedge e^{\sigma}$
in terms of $\Psi_{\ b}^{a}={d}\psi_{\ b}^{a}+\psi_{\ c}^{a}\wedge\psi_{\ b}^{c}$,
where $\mathcal{R}_{\rho\sigma\nu}^{\ \ \ \ \mu}$ is the Riemann
tensor of $\mathcal{Q}$. Note that we use the notation $\psi_{\ b}^{a}=\psi_{\ bc}^{a}e^{c}$
\begin{eqnarray*}
\hat{\Psi}_{\ b}^{a} & = & {d}\hat{\psi}_{\ b}^{a}+\hat{\psi}_{\ c}^{a}\wedge\hat{\psi}_{\ b}^{c}+\hat{\psi}_{\ t}^{a}\wedge\hat{\psi}_{\ b}^{t}\\
 & = & {d}\psi_{\ b}^{a}+{\sqrt{\Lambda}{d}}(\Omega_{\ b}^{a}e^{t})+\psi_{\ c}^{a}\wedge\psi_{\ b}^{c}+\sqrt{\Lambda}\Omega_{\ c}^{a}e^{t}\wedge\psi_{\ b}^{c}+\sqrt{\Lambda}\Omega_{\ b}^{c}\psi_{\ c}^{a}\wedge e^{t}+\Lambda\Omega_{\ [c}^{a}\Omega_{|b|d]}e^{c}\wedge e^{d}\\
 & = & \Psi_{\ b}^{a}+\sqrt{\Lambda}E_{c}(\Omega_{\ b}^{a})e^{c}\wedge e^{t}+\Lambda(\Omega_{\ b}^{a}\Omega_{cd}+\Omega_{\ [c}^{a}\Omega_{|b|d]})e^{c}\wedge e^{d}+\sqrt{\Lambda}(\Omega_{\ b}^{c}\psi_{\ cd}^{a}-\Omega_{\ c}^{a}\psi_{\ bd}^{c})e^{d}\wedge e^{t}\\
 & = & \Psi_{\ b}^{a}+\sqrt{\Lambda}\nabla_{c}\Omega_{\ b}^{a}e^{c}\wedge e^{t}+\Lambda(\Omega_{\ b}^{a}\Omega_{cd}+\Omega_{\ [c}^{a}\Omega_{|b|d]})e^{c}\wedge e^{d}.\\
\hat{\Psi}_{\ a}^{t} & = & {d}\hat{\psi}_{\ a}^{t}+\hat{\psi}_{\ b}^{t}\wedge\hat{\psi}_{\ a}^{b}\\
 & = & \sqrt{\Lambda}E_{[c}(\Omega_{|a|b]})\theta^{c}\wedge\theta^{b}-\sqrt{\Lambda}\Omega_{ab}\psi_{\ c}^{b}\wedge e^{c}+\sqrt{\Lambda}\Omega_{bc}e^{c}\wedge(\psi_{\ a}^{b}+\sqrt{\Lambda}\Omega_{\ a}^{b}e^{t})\\
 & = & \sqrt{\Lambda}(E_{[d}(\Omega_{|a|b]})-\Omega_{ac}\psi_{\ [bd]}^{c}+\Omega_{c[d}\psi_{\ |a|b]}^{c})e^{d}\wedge e^{b}+\Lambda\Omega_{bc}\Omega_{\ a}^{b}e^{c}\wedge e^{t}\\
 & = & \sqrt{\Lambda}\nabla_{[c}\Omega_{|a|d]}e^{c}\wedge e^{d}+\Lambda\Omega_{bc}\Omega_{\ a}^{b}e^{c}\wedge e^{t}.
\end{eqnarray*}
Hence we have that
\begin{eqnarray*}
\mathcal{R}_{cdb}^{\ \ \ a} & = & R_{cdb}^{\ \ \ a}+2\Lambda(\Omega_{\ b}^{a}\Omega_{cd}+\Omega_{\ [c}^{a}\Omega_{|b|d]})\\
\mathcal{R}_{ctb}^{\ \ \ a} & = & \sqrt{\Lambda}\nabla_{c}\Omega_{\ b}^{a}\\
\mathcal{R}_{cda}^{\ \ \ t} & = & 2\sqrt{\Lambda}\nabla_{[c}\Omega_{|a|d]}\\
\mathcal{R}_{cta}^{\ \ \ t} & = & \Lambda\Omega_{bc}\Omega_{\ a}^{b},
\end{eqnarray*}
and thus, using $\mathcal{R}_{\mu\nu}=\mathcal{R}_{\rho\mu\nu}^{\ \ \ \ \rho}$,
\begin{eqnarray*}
\mathcal{R}_{tt} & = & \Lambda\Omega_{bc}\Omega^{bc}=-2n\Lambda=2n\Lambda\mathcal{G}_{tt}\\
\mathcal{R}_{bt} & = & \sqrt{\Lambda}\nabla_{c}\Omega_{\ b}^{c}=0\\
\mathcal{R}_{db} & = & R_{db}+2\Lambda(\Omega_{\ b}^{c}\Omega_{cd}+\frac{1}{2}\Omega_{\ c}^{c}\Omega_{bd}-\frac{1}{2}\Omega_{\ d}^{c}\Omega_{bc})-\Lambda\Omega_{cd}\Omega_{\ b}^{c}\\
 & = & R_{db}+2\Lambda\Omega_{b}^{\ c}\Omega_{dc}\\
 & = & 2(n+1)\Lambda g_{db}-2\Lambda g_{db}=2n\Lambda g_{db}=2n\Lambda\mathcal{G}_{db}.
\end{eqnarray*}

Note that we have used the facts that $g$ is Einstein with Ricci
scalar $4n(n+1)\Lambda$ and that the symplectic form $\Omega$ is divergence-free;
these are justified in the appendix. Since $\mathcal{G}_{at}=0$,
we conclude that
\[
\mathcal{R}_{\mu\nu}=2n\Lambda\mathcal{G}_{\mu\nu}=\frac{\mathcal{R}}{2n+1}\mathcal{G}_{\mu\nu},
\]
i.e. $\mathcal{G}$ is Einstein with Ricci scalar $2n(2n+1)\Lambda$.

\begin{flushright}
$\square$
\par\end{flushright}

Physically, this is a Kaluza-Klein reduction with constant dilation
field and where the Maxwell two-form is related to the reduced metric
by $\Omega_{a}^{\ c}\Omega_{cb}=g_{ab}$. This is what allows both
the reduced and lifted metric to be Einstein. A more general discussion
can be found in \cite{Pope}.

From the Cartan perspective, $\mathcal{G}_{\Lambda=1}$ can be thought
of as a metric on the $2n+1$-dimensional space obtained by taking
a quotient $\tilde{q}:P\mapsto P/SL(n,\mathbb{R})$ of the Cartan
bundle, where we embed $SL(n,\mathbb{R})\subset GL(n,\mathbb{R})$
in $H$ as in (\ref{eq:gl(n)_embed_H}) but with $a$ now denoting
an element of $SL(n,\mathbb{R})$ (so that $\mathrm{det}(a)^{-1}=1$).
This new subgroup acts adjointly on $\theta$ as
\[
\begin{pmatrix}1 & 0\\
0 & a
\end{pmatrix}\begin{pmatrix}-\mathrm{tr}(\phi) & \eta\\
\omega & \phi
\end{pmatrix}\begin{pmatrix}1 & 0\\
0 & a^{-1}
\end{pmatrix}=\begin{pmatrix}-\mathrm{tr}(\phi) & \eta a^{-1}\\
a\omega & \phi
\end{pmatrix},
\]
so not only is the inner product $\eta\omega$ invariant but also
the $(0,0)$-component $\theta_{\ 0}^{0}=\mathrm{-tr}\phi$, which
is a scalar one-form whose exterior derivative is constrained by (\ref{eq:curvature_2-form})
to be $\mathrm{d}\theta_{\ 0}^{0}=-\theta_{\ i}^{0}\wedge\theta_{\ 0}^{i}=-\mathrm{Ant}(\eta\wedge\omega)$.
Thus, denoting by ${A}$ the object on $\mathcal{Q}=P/SL(n,\mathbb{R})$
which is such that $\tilde{q}^{*}A=\mathrm{tr}\phi$, we have that
${d}A=\Omega$ (where we are now taking $\Omega$ and $g$
to be defined on $\mathcal{Q}$ by $\tilde{q}^{*}\Omega=\mathrm{Ant}(\eta\wedge\omega)$
and $\tilde{q}^{*}g=\mathrm{Sym}(\eta\wedge\omega)$ respectively,
or equivalently redefining $\tilde{\Omega}=\sigma^{*}\Omega$ and
$\tilde{g}=\sigma^{*}g$).

We then have a natural way of constructing a metric $\mathcal{G}$
on $\mathcal{Q}$ as a linear combination of $g$ and $e^{t}\odot e^{t}$,
where $e^{t}$ is $A$ up to addition of some exact one-form. It turns
out that the choice of linear combination such that $\mathcal{G}$
is Einstein is 
\[
\mathcal{G}=g-e^{t}\odot e^{t}.
\]
The fact that this metric is exactly (\ref{eq:G}) can be verified
by constructing the Cartan connection of $(S,[\nabla])$ explicitly
in terms of a representative connection $\nabla\in[\nabla]$.

\section*{Appendix A: Ricci scalar of $(M,g)$ and divergence of $\Omega$}
\appendix
We calculate these using the Cartan formalism, again using the basis
(\ref{eq:basis}). In this basis we have $g$ as above (\ref{eq:g_cov_const})
and
\[
\Omega=\sum_{A=0}^{n-1}e^{A}\wedge e^{A+n}\quad\implies\quad\Omega_{ab}=\sum_{A=0}^{n-1}\delta_{[a}^{A}\delta_{b]}^{A+n}.
\]
Note that from now on we will omit the summation sign and use the
summation convention regardless of whether $A, B$--indices are up or
down. As in \S \ref{kk_curved_case}, we look for $\psi_{\ b}^{a}$
by considering ${d}e^{a}$ (recall that $A,B=0,\dots,n-1$ and
$a,b=0,\dots,2n-1$):
\begin{eqnarray*}
{d}e^{A} & = & 0\\
{d}e^{A+n} & = & -(E_{D}(\Gamma_{AB}^{C})z_{C}-\Lambda^{-1}E_{D}(P_{AB})){d}x^{D}\wedge{d}x^{B}-(\Gamma_{AB}^{C}-2\Lambda z_{(A}\delta_{B)}^{C}){d}z_{C}\wedge{d}x^{B}\\
 & = & -(E_{D}(\Gamma_{AB}^{C})z_{C}-\Lambda^{-1}E_{D}(P_{AB}))e^{D}\wedge e^{B}\\
 &  & -(\Gamma_{AB}^{C}-2\Lambda z_{(A}\delta_{B)}^{C})(e^{C+n}+(\Gamma_{CE}^{D}z_{D}-\Lambda z_{C}z_{E}-\Lambda^{-1}P_{CE})e^{E})\wedge e^{B}\\
 & = & \bigl[\Lambda^{-1}E_{E}(P_{AB})-E_{E}(\Gamma_{AB}^{C})z_{C}+\Lambda^{-1}\Gamma_{AB}^{C}P_{CE}-\Gamma_{AB}^{C}\Gamma_{CE}^{D}z_{D}+\Lambda\Gamma_{AB}^{C}z_{E}z_{C}\\
 &  & +2\Lambda z_{(A}(\Gamma_{B)E}^{D}z_{D}-\Lambda z_{B)}z_{E}-\Lambda^{-1}P_{B)E})\bigr]e^{E}\wedge e^{B}+(2\Lambda z_{(A}\delta_{B)}^{C}-\Gamma_{AB}^{C})e^{C+n}\wedge e^{B}\\
 & = & \bigl[\Lambda^{-1}D_{E}P_{AB}-(D_{E}\Gamma_{AB}^{C})z_{C}-2z_{(A}P_{B)E}\bigr]e^{E}\wedge e^{B}+(2\Lambda z_{(A}\delta_{B)}^{C}-\Gamma_{AB}^{C})e^{C+n}\wedge e^{B}\\
\end{eqnarray*}
Note that we have used $D$ to denote the chosen connection on $S$
with components $\Gamma_{BC}^{A}$. Next we wish to read off the spin
connection $\psi_{\ b}^{a}$ such that ${d}e^{a}=-\psi_{\ b}^{a}\wedge e^{b}$
and the following index symmetries are satisfied:
\begin{eqnarray*}
\psi_{\ B}^{A} & = & \frac{1}{2}\psi_{A+n\, B}=-\frac{1}{2}\psi_{B\, A+n}=-\psi_{\ A+n}^{B+n}\\
\psi_{\ B+n}^{A} & = & \frac{1}{2}\psi_{A+n\, B+n}=-\frac{1}{2}\psi_{B+n\, A+n}=-\psi_{\ A+n}^{B}\\
\psi_{\ B}^{A+n} & = & \frac{1}{2}\psi_{A\, B}=-\frac{1}{2}\psi_{B\, A}=-\psi_{\ A}^{B+n}
\end{eqnarray*}
We find that 
\begin{eqnarray*}
\psi_{\ C+n}^{A+n} & = & (2\Lambda z_{(A}\delta_{B)}^{C}-\Gamma_{AB}^{C})e^{B}=-\psi_{\ A}^{C}\\
\psi_{\ B}^{A+n} & = & \bigl[2(D_{[A}\Gamma_{B]C}^{D})z_{D}-2\Lambda^{-1}D_{[A}P_{B]C}^{\mathrm{S}}-\Lambda^{-1}D_{C}P_{AB}^{\mathrm{A}}+2z_{(B}P_{C)A}-2z_{(A}P_{C)B}\bigr]e^{C}=:A_{ABC}e^{C}\\
\psi_{\ B+n}^{A} & = & 0.
\end{eqnarray*}
One can check that these satisfy both the index symmetries above and
are such that ${d}e^{a}=-\psi_{\ b}^{a}\wedge e^{b}$, and
we know from theory that there is a unique set of $\psi_{\ b}^{a}$
that have both of these properties. Note that we have used $P^{\mathrm{S}}$
and $P^{\mathrm{A}}$ to denote the symmetric and antisymmetric parts
of $P$ in order to avoid too much confusion from having multiple
symmetrisation brackets in the indices.

We are now ready to calculate the divergence of $\Omega$. Since it
is also covariantly constant in this basis, we obtain
\[
\nabla_{c}\Omega_{ab}=-\psi_{\ ac}^{d}\Omega_{db}-\psi_{\ bc}^{d}\Omega_{ad}=-\psi_{\ ac}^{d}\Omega_{db}+\psi_{\ bc}^{d}\Omega_{da}=2\Omega_{d[a}\psi_{\ b]c}^{d}.
\]
We can split the right hand side into
\begin{eqnarray*}
\Omega_{da}\psi_{\ bc}^{d} & = & \Omega_{Ca}\psi_{\ bc}^{C}+\Omega_{C+n\, a}\psi_{\ bc}^{C+n}\\
 & = & \delta_{[C}^{A}\delta_{a]}^{A+n}\psi_{\ bc}^{C}+\delta_{[C+n}^{A}\delta_{a]}^{A+n}\psi_{\ bc}^{C+n}\\
 & = & \frac{1}{2}\Bigl(-\delta_{a}^{C+n}\delta_{b}^{A}\delta_{c}^{B}(2\Lambda z_{(A}\delta_{B)}^{C}-\Gamma_{AB}^{C})-\delta_{a}^{C}\delta_{b}^{D+n}\delta_{c}^{B}(2\Lambda z_{(C}\delta_{B)}^{D}-\Gamma_{CB}^{D})-\delta_{a}^{C}\delta_{b}^{D}\delta_{c}^{E}A_{CDE}\Bigr).\\
\end{eqnarray*}
The first two terms are the same but with $a\leftrightarrow b$, so
are lost in the antisymmetrisation. Thus
\[
\nabla_{c}\Omega_{ab}=-\delta_{[a}^{C}\delta_{b]}^{D}\delta_{c}^{E}A_{CDE}.
\]
Tracing amounts to contracting this with $g^{ac}$:
\[
\nabla^{c}\Omega_{cb}=-\delta_{[a}^{C}\delta_{b]}^{D}g^{ac}\delta_{c}^{E}A_{CDE}=-\delta_{[a}^{C}\delta_{b]}^{D}g^{aE}A_{CDE},
\]
but $g^{aE}$ is non-zero only when $a=E+n>n$ and $\delta_{[a}^{C}\delta_{b]}^{D}$
is non-zero only when $a=C\leq n$ or $a=D\leq n$. We can therefore
conclude that the right hand side is zero and $\Omega$ is divergence-free.

Finally, we calculate the Ricci scalar of $g$ (given that it's Einstein)
via the curvature two-forms $\Psi_{\ b}^{a}={d}\psi_{\ b}^{a}+\psi_{\ c}^{a}\wedge\psi_{\ b}^{c}=\frac{1}{2}R_{cdb}^{\ \ \ a}e^{c}\wedge e^{d}$.
We are only interested in non-zero components of the Ricci tensor
such as $R_{A\, B+n}=R_{cA\, B+n}^{\ \ \ \ \ \ \ \ c}$. In fact, we will
calculate only $R_{E+n\, B}$, for which we need to consider $R_{D\, E+n\, B}^{\ \ \ \ \ \ \ \ \ A}$
and $R_{D+n\, E+n\, B}^{\ \ \ \ \ \ \ \ \ \ \ \ A+n}$, i.e. we need only calculate
$\Psi_{\ B}^{A}$ and $\Psi_{\ \ B}^{A+n}$. 
\[
\Psi_{\ B}^{A}={d}\Bigl((\Gamma_{BC}^{A}-2\Lambda z_{(B}\delta_{C)}^{A})e^{C}\Bigr)+\psi_{\ C}^{A}\wedge\psi_{\ B}^{C}+\psi_{\ C+n}^{A}\wedge\psi_{\ B}^{C+n}.
\]
The last term vanishes since $\psi_{\ B}^{C+n}=0$, and the middle
term only has components that look like $\frac{1}{2}R_{DEB}^{\ \ \ \ \ A}e^{D}\wedge e^{E}$,
so the only term we are interested in is 
\[
-2\Lambda{d}z_{(B}\delta_{C)}^{A}e^{C}=-2\Lambda\delta_{(C}^{A}(e^{B)+n}+(\Gamma_{B)E}^{D}z_{D}-\Lambda z_{B)}z_{E}-\Lambda^{-1}P_{B)E})e^{E})\wedge e^{C}.
\]
Again, discarding the $e^{E}\wedge e^{C}$ term gives
\[
-\Lambda(e^{B+n}\wedge e^{C}+\delta_{B}^{A}e^{C+n}\wedge e^{C})=\frac{1}{2}R_{D\, E+n\, B}^{\ \ \ \ \ \ \ \ \ A}e^{D}\wedge e^{E+n}+\frac{1}{2}R_{E+n\, D\, B}^{\ \ \ \ \ \ \ \ \ A}e^{E+n}\wedge e^{D},
\]
so we conclude
\[
R_{D\, E+n\, B}^{\ \ \ \ \ \ \ \ \ A}=\Lambda(\delta_{B}^{A}\delta_{D}^{E}+\delta_{D}^{A}\delta_{B}^{E}).
\]
The other Riemann tensor component we need to know to calculate $R_{E+n\, B}=R_{c\, E+n\, B}^{\ \ \ \ \ \ \ \ \ c}$
is $R_{D+n\, E+n\, B}^{\ \ \ \ \ \ \ \ \ \ \ \ A+n}$, so we examine
\[
\Psi_{\ \ B}^{A+n}={d}\psi_{\ \ B}^{A+n}+\psi_{\ \ C}^{A+n}\wedge\psi_{\ B}^{C}+\psi_{\ \ C+n}^{A+n}\wedge\psi_{\ \ B}^{C+n},
\]
but none of these terms have $e^{D+n}\wedge e^{E+n}$ components,
so $R_{D+n\, E+n\, B}^{\ \ \ \ \ \ \ \ \ \ \ \ A+n}=0$. Hence 
\[
R_{E+n\, B}=\delta_{D}^{A}R_{D\, E+n\, B}^{\ \ \ \ \ \ \ \ \  A}=\Lambda(\delta_{B}^{E}+n\delta_{B}^{E})=\Lambda(n+1)\delta_{B}^{E}.
\]
Setting this equal to $\frac{R}{2n}g_{E+n\, B}=\frac{R}{4n}\delta_{B}^{E}$
we find 
\[
R=4n(n+1)\Lambda,
\]
as required.
\section*{Appendix B: $SU(\infty)$ Toda from an ALH instanton}
\appendix

Consider the  Hyper--K\"ahler metric in the 
Gibbons--Hawking class \cite{Gibbons:1979zt}, where the harmonic potential 
on $\R^3$ is a linear function, i. e.
\be
\label{GH}
g = z \left( dx^2+dy^2+dz^2 \right) + z^{-1} \Big( dt + 
\frac{1}{2}(xdy-ydx) \Big)^2.
\ee
It admits a 4-dimensional group of isometries generated by the Killing 
vectors
\[
\p_t,\quad x\p_y-y\p_x,\quad \p_x-\frac{y}{2}\p_t, \quad
\p_y+\frac{x}{2}\p_t.
\]
Introducing the radial coordinate $r$ by $z=(9/4)^{1/3}r^{2/3}$
we find the volume growth $r^{4/3}$, so the metric is ALH in the sense of 
\cite{ALHref}
(although it is not complete, as it is singular at $r=0$). For the ALH behaviour we need to identify
$(x, y)$ with coordinates on a two-torus (see \cite{Gibbons:1998ie} for 
another interpretation of this solutions in terms of BPS domani walls).
\subsection*{Reduction by $K=x\p_y-y\p_x$.} This Killing vector is not compatible with the toric topology of the $(x-y)$ space, but it will lead to a 
non--trivial solution to (\ref{md_toda}) with $\epsilon=-1$.

Set $x=(2/3)^{1/3}\rho\cos{\theta}, y=(2/3)^{1/3}\rho\sin{\theta},
\tau=(3/2)^{1/3}X$. The K\"ahler form is
\[
\omega=r^{-1/3}dr\wedge(dX+\frac{1}{3}\rho^2d\theta)+r^{2/3}\rho d\rho\wedge 
d\theta.
\]
It is preserved by the Killing vector $K=\p/\p\theta$. Formula (\ref{ztilde})
yields the moment map
\[
Z=\frac{1}{2}r^{2/3}\rho^2.
\]
We eliminate $\rho$ in terms of $Z$ and $r$ and define
\[
Y=\frac{1}{3}\frac{Z}{r^{2/3}}-\frac{3}{4}r^{4/3}.
\]
The resulting solution of the Toda equation \ref{md_toda} with $\epsilon=-1$ is given by
$e^U=2r^{-2/3}Z$ . It is 
constant on a parabolic cylinder
$
(e^U/6-Y)e^{2U}=3Z^2.
$
It is given by
\[
e^U=
\sqrt [3]{8\,{Y}^{3}+9\,{Z}^{2}+3\,\sqrt {16\,{Z}^{2}{Y}^{3}+9\,{Z}^{4
}}}+{\frac {4{Y}^{2}}{\sqrt [3]{8\,{Y}^{3}+9\,{Z}^{2}+3\,\sqrt {16\,
{Z}^{2}{Y}^{3}+9\,{Z}^{4}}}}}+2\,Y.
\]
(See \cite{prim} for other solutions constant on cylinders).
\subsection*{Reduction by $K=\p_x+\p_y+\frac{1}{2}(x-y)\p_t$}
This  Killing vector is
compatible with the toric structure.
Following the procedure above, then one finds the solution to 
(\ref{md_toda}) to be $u=0$. This is because $K$ 
is a tri--holomorphic Killing vector.

\end{document}